\documentclass[a4paper,12pt]{amsart}

\usepackage{amsmath, amsfonts, amsthm,amssymb}
\usepackage[francais]{babel}
\usepackage[latin1]{inputenc}

\newcommand{\Ker}{\textrm{Ker}}
\newcommand{\disc}{\textrm{disc}}
\newcommand{\Fred}{\textrm{Fred}}
\newcommand{\rr}{\textrm{red}}
\newcommand{\tr}{\textrm{tr}}
\newcommand{\red}{\textrm{Red}}
\newcommand{\res}{\textrm{res}}
\newcommand{\cond}{\textrm{cond}}

\newcommand{\End}{\textrm{End}}
\newcommand{\Comp}{\textrm{Comp}}

\newcommand{\FF}{F}

\newcommand{\AAA}{\mathbb{A}}

\newcommand{\TTT}{\mathbb{M}} 
\newcommand{\OO}{\mathcal{O}}

\newcommand{\JL}{\textrm{JL}}
\newcommand{\Tr}{\textrm{Tr}}
\newcommand{\DE}{D^{0,\varepsilon,d}}
\newcommand{\DQ}{D^{D,\varepsilon}}
\newcommand{\WW}{\mathcal{W}}

\newcommand{\NN}{\mathbb{N}}

\newcommand{\ZZ}{\mathbb{Z}}
\newcommand{\QQ}{\mathbb{Q}}
\newcommand{\RR}{\mathbb{R}}
\newcommand{\CC}{\mathbb{C}}
\newcommand{\C}{\mathcal{C}}
\newcommand{\HH}{\mathcal{H}}
\newcommand{\PP}{\mathbb{P}}

\newcommand{\EEE}{\mathbb{E}}
\newcommand{\divi}{\, |\!| \,}
\newcommand{\diag}{\textrm{diag}}
\newcommand{\ME}{F^{0,\varepsilon,d}}
\newcommand{\MQ}{F^{D,\varepsilon}}

\newtheorem{thm}[subsection]{Théorème}
\newtheorem{prop}[subsection]{Proposition}
\newtheorem{cor}[subsection]{Corollaire}
\newtheorem{lemme}[subsection]{Lemme}

\def\]{\textup{\mbox{]\hspace{-.15em}]}}}
\def\[{\textup{\mbox{[\hspace{-.15em}[}}}

\input xy
\xyoption{all}

\title{Une correspondance de Jacquet-Langlands $p$-adique}

\begin{document}

\maketitle

\begin{center} Gaëtan CHENEVIER \end{center}

\vspace{3 mm}

\begin{center} Ecole normale supérieure, DMA \footnote{chenevie@dma.ens.fr} \end{center}

\begin{center} Université Paris 7 \end{center}

\vspace{1 cm}

{\it Abstract:} In this paper, we extend the Jacquet-Langlands'correspondence, between Hecke-modules of usual modular forms and quaternionic modular
forms, to overconvergent $p$-adic forms of finite slope. We show that this
correspondence respects $p$-adic families and is induced by
an isomorphism between some associated eigencurves. \par \vspace{5 mm}

\par \vspace{2 mm}

\section{Introduction}

	Soit $D/\QQ$ l'algèbre de quaternions définie de discriminant $d$,
$p$ un nombre premier et $N$ un entier tels que $(Np,d)=(N,p)=1$, nous
établissons dans ce texte un transfert "à la Jacquet-Langlands" bijectif
entre les formes modulaires $p$-adiques surconvergentes de pente finie (\cite{K},
\cite{Col2}), propres, cuspidales et nouvelles en $d$, et les formes modulaires
$p$-adiques quaternioniques pour $D$ de pente finie (\cite{Buz2}), en niveau
modéré $N$ et poids-caractère quelconque (théorème \ref{jlp}). Ce transfert est
Hecke-équivariant et coïncide avec la correspondance de Jacquet-Langlands
usuelle lorsqu'on le restreint aux formes modulaires "classiques" de part et
d'autre. Mieux, il respecte les familles $p$-adiques, et nous prouvons qu'il
provient d'un isomorphisme rigide-analytique entre les courbes de Hecke ("the eigencurves")
correspondantes (théorème \ref{jlpfamille}). \par \vspace{2 mm}
	Notre preuve est basée sur l'existence de systèmes de modules de
Banach de part et d'autre (\ref{banachsystem}), et de leur comparaison. On
démontre en fait un énoncé général (théorème \ref{general}) sur la comparaison de
ces systèmes. Des arguments de Zariski-densité de points "classiques" sur certaines
hypersurfaces de Fredholm et sur les variétés de Hecke y jouent un rôle important. On en déduit notre correspondance en utilisant la
correspondance de Jacquet-Langlands usuelle, et les assertions de "classicité
en petite pente" de part et d'autre. \par \vspace{2 mm}
	Un des charmes de cette correspondance est que les objets qu'elle
compare sont, dans une large mesure, de natures assez différentes. Les formes
modulaires $p$-adiques surconvergentes usuelles sont des sections
de fibrés surconvergentes sur le lieu ordinaire des courbes modulaires, alors que les formes modulaires
quaternioniques $p$-adiques sont purs produits de la théorie des groupes,
spécialement des représentations $p$-adiques de $GL_2(\QQ_p)$. Nous y voyons
plusieurs intérêts, l'un d'entre eux étant la possibilité d'introduire des méthodes de
théorie des représentations $p$-adiques des groupes de Lie $p$-adiques dans la théorie des formes modulaires surconvergentes. 
Ceci est par exemple en accord avec une philosphie de Langlands purement $p$-adique
(encore non formulée à notre connaissance). Un autre intêret vient de qu'il reste une multitude de questions non résolues concernant "the
eigencurve" (voir l'introduction de \cite{eigen}), l'origine plus combinatoire des formes quaternioniques peut
permettre, sinon de les résoudre, de faciliter tout au moins les expériences numériques.
Nous discutons de certaines conséquences de notre correspondance, ainsi que
de problèmes encore ouverts, dans la dernière partie du texte. \par
\vspace{2 mm}
Dans \cite{ch}, des variétés de Hecke sont construites pour tous les groupes algébriques $G$ sur
$\QQ$ tels que $G(\RR)=U_n(\CC)$, $G(\QQ_p)=GL_n(\QQ_p)$, et on peut se demander
de manière générale si les transferts usuels se prolongent, et comment, aux formes modulaires surconvergentes
pour ces groupes. En ce qui concerne les correspondances de Jacquet-Langlands
éventuelles entre de tels groupes, nos résultats s'y appliquent sans modification, 
du moment que la correspondance classique entre les deux groupes en question est connue 
pour suffisament de poids (voir le théorème \ref{general}). \par \vspace{2 mm}
	Ce travail repose substantiellement sur ceux de Buzzard, Coleman,
et Mazur, nous les en remercions. La question de l'existence ou non d'un morphisme entre les deux courbes de Hecke avait notamment été posée par K.Buzzard.

\par \vspace{3 mm}

{\bf Notations}: $p$ est un nombre premier. 
$\CC_p$ est le complété d'une clôture algébrique de $\QQ_p$, muni
de la norme telle que $|p|=1/p$, $v(.)$ est la valuation associée. Si $X/\CC_p$ est un espace rigide, $A(X)$ est
l'anneau des fonctions rigides analytiques sur $X$; si $X$ est de plus
affinoide réduit, la norme sup. sur $X$ munit $A(X)$ d'une structure de $\CC_p$-algèbre de Banach.
$\AAA^1$ désigne la droite rigide sur $\CC_p$. $\AAA$ (resp. $\AAA_f$, resp. $\AAA_f^{(d)}$) est l'anneau des adèles (resp.
adèles finies, resp. adèles hors de $\{d,\infty\}$) de $\QQ$. Si $d$ et $n$
sont des entiers, on note $d\divi n$ si $d$ divise $n$ et $(d,n/d)=1$. Si $B$ est un anneau, $B^*$ désigne le groupe multiplicatif de ses
inversibles. $\diag(a,b)$ désigne la matrice $\left( \begin{array}{cc} a &  0 \\  0 & b \end{array}
\right)$. 

\tableofcontents

\section{Systèmes de modules de Banach} 

\subsection{Modules de Banach} \label{banachsystem} La notion de systèmes de modules de Banach
orthonormalisables est introduite dans \cite{eigen} \S 4.3. \par \vspace{1 mm}
On appellera {\it système d'espaces de Banach} 
la donnée d'un ensemble $E=\{ E_n, i_n, n \in \NN\}$ de $\CC_p$-espaces de Banach
orthonormalisables $E_n$, et d'applications $\CC_p$-linéaires
compactes $i_{n}: E_n \rightarrow E_{n+1}$, on notera $E^{\dagger}$ la
limite inductive des $E_n$ selon les $i_n$. Si $\WW/\CC_p$ est un espace rigide réduit, on appellera {\it faisceau de modules de Banach
orthonormalisables sur $\WW$} la donnée d'un faisceau $B$ sur $\WW$ tel que
pour tout ouvert affinoide $V \subset \WW$, $B(V)$ ait une structure de
$A(V)$-module de Banach orthonormalisable, et tel que si $V \subset V'$ sont
des ouverts affinoides de $\WW$, l'application canonique
$B(V')\widehat{\otimes}_{A(V')}A(V) \rightarrow B(V)$ soit un isomorphisme
de $A(V)$-modules de Banach. Un {\it système de modules de Banach sur $\WW$}
est la donnée d'un ensemble $E=\{E_n,i_n, n\in \NN\}$ où $E_n$ est un faisceau de modules de
Banach orthonormalisables sur $\WW$, et $i_n: E_n \rightarrow E_{n+1}$
est un morphisme de faisceaux tel que si $V$ est un ouvert affinoide de
$\WW$, $i_n(V): E_n(V) \rightarrow E_{n+1}(V)$ soit $A(V)$-linéaire
compacte. On notera alors $E(V,n):=E_n(V)$, et on abrégera le tout en
$E=(E(V,n))$. Enfin, on dispose d'une notion évidente de sous-système de
modules de Banach sur $\WW$. \par \vspace{2 mm}

{\it Exemple:} Si $E=\{ E_n, i_n, n \in \NN\}$ est un système d'espaces de
Banach, on peut lui associer un système de modules de Banach sur $\WW$ comme
suit. Soit $E_{\WW,n}$ le faisceau de modules de Banach
orthonormalisables associé
à $E_n$ au sens de \cite{ch} 3.7.2. On rappelle que si $V$ est ouvert affinoide, on a
$E_{\WW,n}(V):=A(V)\widehat{\otimes}_{\CC_p} E_n$ que l'on notera encore
$E(V,n)$. $i_n$ nous fournit de plus par extension des scalairs une application compacte $A(V)$-linéaire
$i_n(V): E(V,n) \rightarrow E(V,n+1)$. On
appellera $(E(V,n))$ {\it le système de modules de Banach sur $\WW$ associé à
$E$}. \par \vspace{2 mm}

Fixons $E=(E(V,n))$ un système de modules de Banach sur $\WW$. 
Si $x \in \WW(\CC_p)$, on notera $E_{x}$ le système d'espaces de
Banach déduit de $E$ par évaluation en $x$:
$E_{x,n}:=E(V,n)\widehat{\otimes}_{A(V)}\CC_p$ pour tout $V$ tel que $x \in
V(\CC_p)$ ($A(V) \rightarrow \CC_p$ étant l'évaluation en $x$), $i_{x,n}: E_{x,n}
\rightarrow E_{x,n+1}$ l'application compacte déduite. Dans nos applications, les $i_{x,n}$ seront toujours injectives. C'est
par exemple le cas si $E_{\WW}$ est le système de modules de Banach sur $\WW$ associé
a un système d'espaces de Banach $\{E_n, \, i_n, \, n \in \NN\}$ tel que les
$i_n$ soient injectives. \par \vspace{2 mm}

	On appellera endomorphisme de $E$ la donnée, pour chaque $V \subset \WW$
ouvert affinoide, et pour tout $n$ assez grand ($V$ étant fixé), d'un endomorphisme continu
$A(V)$-linéaire $U(V,n)$ de $E(V,n)$, tel que les $U(V,n)$ commutent aux
applications $E(V,n) \rightarrow E(V,n+1)$, et aux changements de bases
ouverts affinoides $E(V,n) \rightarrow E(V',n)$ lorsqu'ils sont définis.
On dira que $U:=(U(V,n))$ est un endomorphisme de $E$ et on identifiera deux
endomorphismes $U$ et $U'$ si pour tout $V$ et $n$ assez grand,
$U(V,n)=U'(V,n)$. L'ensemble $\End(E)$ de ces endomorphismes est alors une 
$\CC_p$-algèbre de manière naturelle. Si $H$ est un anneau (resp. $G$ un
moinoide), une représentation de
$H$ (resp. de $G$) sur $E$ est la donnée d'un morphisme de $\CC_p$-algèbre $H \rightarrow
\End(E)$ (resp. $\CC_p[G] \rightarrow \End(E)$). On parlera alors de système
de $H$-modules de Banach sur $\WW$. \par \vspace{2 mm}

	On note $\Comp(E)$ l'idéal bilatère de $\End(E)$ composé des
éléments ayant la propriété suivante: pour $V \subset \WW$ un ouvert
affinoide fixé, il existe pour tout $n$ assez grand un endomorphisme
$A(V)$-linéaire continu $T(V,n): E(V,n+1) \rightarrow E(V,n)$ tel que 
le diagramme suivant soit commutatif:

$$\xymatrix{ E(V,n) \ar@{->}[d]_{i_{n}(V)} \ar@{->}[rr]^{U(V,n)}
& & E(V,n) \ar@{->}[d]^{i_{n}(V)} \\
 E(V,n+1) \ar@{->}[rr]_{U(V,n+1)} \ar@{->}[urr]^{T(V,n)} & & E(V,n+1) }$$

En particulier, $U(V,n)$ est compact, et \cite{BMF} A.2.3
assure que $\det(1-U(V,n)_{|E(V,n)}) \in 1+TA(V)\{\{T\}\}$ est indépendant de $n$ assez
grand. La formation des séries de Fredholm commutant aux changements de base
ouverts affinoides, toutes ces séries proviennent par restriction d'une
unique série de Fredholm \Fred$_E$($U$) sur $\WW$, \Fred$_E$($U$) $\in
1+TA(\WW)\{\{T\}\}$. \par \vspace{2 mm}

\subsection{L'espace des poids-caractères} \label{poids} ${}^{}$ \par
\vspace{1 mm}
\par On fixe $p>2$ dans ce qui suit, $\Lambda$ est l'anneau local complet $\ZZ_p[[(1+p\ZZ_p)]]$,
$\WW$ la boule ouverte rigide de centre $1$ et rayon $1$, $\WW(\CC_p):=\{z \in \CC_p, \, \, |z-1|<1\}$. L'application $\Lambda \rightarrow
A(\WW)$ définie par $[1+p] \mapsto Z$ identifie 
topologiquement $\Lambda$ aux fonctions analytiques $\QQ_p$-valuées sur $\WW(\QQ_p)$ bornées par
$1$ sur tout $\WW(\CC_p)$. L'application $\kappa \mapsto \kappa(1+p)$ induit une
bijection:
$$\textrm{Hom}_{gr-cont}(1+p\ZZ_p,\CC_p^*) \simeq \WW(\CC_p)$$  

On dispose d'un caractère "universel" continu déterminé par
$\kappa^{univ}(1+p)=Z$: $$\kappa^{univ}: 1+p\ZZ_p
\rightarrow A(\WW)^*$$

On note $\mu_{p^{\infty}}:=\{\zeta \in \CC_p^*, \, \, \exists r \in \NN, 
\zeta^{p^r}=1\}$, $\mu_{p^{\infty}} \subset \WW(\CC_p)$. Si $\zeta
\in \CC_p^*$ est tel que $\zeta^{p^r}=1$,  $\zeta.(1+p)^k \in \WW(\CC_p)$ correspond 
au caractère $x \rightarrow x^k \chi(x)$, $\chi$
étant le caractère d'ordre fini de $1+p\ZZ_p$ tel que $\chi(1+p)=\zeta$,
caractère que l'on notera $(k,\chi)$. Dans ce cas, on appellera conducteur de
$\kappa$, noté $\cond(\chi)$, le plus petit
entier naturel $r$ tel que $\chi$ soit trivial sur $1+p^r\ZZ_p$ . \par
\vspace{1 mm}
On aura besoin en \ref{quat} d'introduire pour $r \geq 1$, la boule rigide
$\WW_r \subset \WW$: 
 $$\WW_r(\CC_p):=\{ \kappa \in \WW(\CC_p), \, \,
|\kappa-1|< p^{-\frac{1}{p^{r-1}(p-1)}} \}$$ On a par exemple $(1+p)^k\zeta \in
\WW_r$
quand $\zeta^{p^{r-1}}=1$. Tout $\kappa$ dans $\WW_r(\CC_p)$ est un caractère de
$1+p\ZZ_p$ de restriction analytique à
$1+p^r\ZZ_p$. Mieux, si $V \subset \WW_r$ est ouvert affinoide, la restriction de $\kappa^{univ}$ à
$A(V)$ est un
caractère $A(V)$-valué de $1+p\ZZ_p$ dont la restriction à $1+p^r\ZZ_p$
est analytique. \par \vspace{1 mm}
	On note $\tau: (\ZZ/p\ZZ)^* \rightarrow \ZZ_p^*$ le caractère
de Teichmüller. On verra en général les caractères de $1+p\ZZ_p$ comme des
caractères de $\ZZ_p^*$, en les étendant trivialement sur les racines de
l'unité. 
\par \vspace{3 mm}

\section{Formes modulaires $p$-adiques}

\label{formes}  Nous allons dans ce qui suit rappeler les acteurs essentiels de ce papier. On fixe un nombre premier
$p\geq 5$, des entiers $N$ et $d$, $(N,p)=(N,d)=(p,d)=1$, $d$ sans facteur
carré.  \par \vspace{1 mm}

Soit $\mathcal{H}$ la $\ZZ$-algèbre commutative de polynômes sur les lettres $S_l$,
$T_l$, si $l$ est premier et $(l,Ndp)=1$, et $U_l$ si $l$ premier divise
$Ndp$. On fixe $\varepsilon=\varepsilon_p\varepsilon_N$ un caractère de $(\ZZ/p\ZZ\times
\ZZ/N\ZZ)^*$. On notera aussi, si $n \geq 1$, $T_n$ l'élément de $\HH$ obtenu par les
formules usuelles:

$$\sum_{n\geq 1} T_n n^{-s} =\prod_{l|Npd} (1-U_ll^{-s})^{-1}\prod_{l
\nmid Npd} (1-T_ll^{-s}+l S_ll^{-2s})^{-1}$$
\par \vspace{3 mm}

\subsection{Formes modulaires surconvergentes} ${}^{}$ \par
\vspace{2 mm}
\noindent On fait des rappels sur certaines constructions faites dans
\cite{eigen} \S 2.1,
\cite{BMF}
\S 4.3 . \par \vspace{2 mm}  Soit $X_1(Np,d)$ la
courbe propre et plate sur $\ZZ_p$ paramétrant les courbes elliptiques
généralisées avec structure de niveau de type $\Gamma_1(Np)\cap\Gamma_0(d)$, $\omega$ le
faisceau inversible habituel sur $X_1(Np,d)$. Si $p>3$, on 
rappelle que la série d'Eisenstein normalisée de niveau $1$, $E_{p-1}$, est une section globale de
$\omega^{p-1}$ sur $X_1(Np,d)/\ZZ_p$ qui relève l'invariant de Hasse. 
Si $0\leq v<1 \in \QQ$, on définit $X_1(Np,d)(v)$ comme étant
l'ouvert affinoide de $X_1(Np,d)^{rig}$ sur lequel $|E_{p-1}|\geq p^{-v}$,
et $Z_1(Np,d)(v)$ la composante connexe affinoide de $\infty$ dans
$X_1(Np,d)(v)$ (\cite{eigen} \S 2.1, \S 4.3). Pour tout $m \geq 1$, on s'intéressera plus généralement (\cite{eigen} \S2)
à la courbe rigide analytique (lisse irréductible) $Z_1(Np^m,d)(v)$, avec  $0\leq v <
1$, qui est la composante connexe affinoide de $\infty$ dans $X_1(Np^m,d)(v)$.
$$M_k(p^m,d,v):=H^0(Z_1(Np^m,d)(v),\omega^k)$$ est le $\CC_p$-espace de
Banach des formes modulaires $v$-surconvergentes de poids $k$, de niveau
$\Gamma_1(Np^m)\cap \Gamma_0(d)$, il est muni d'une action naturelle de $\HH$ (\cite{eigen} \S
3.2, \cite{BMF} B5).
On fixe ici, et pour tout le texte, une suite réelle décroissante $(v_n)_{n\in \NN}$, telle que
$\forall n \in \NN, \, v_n=p^{-n}v_0 \in [\frac{p}{p^{n+1}(p+1)},\frac{p}{p^{n}(p+1)}[\cap
\QQ$. La construction de nos systèmes de modules de Banach dépend de
$v_0$ d'une manière qui nous importera peu, pour alléger les notations
nous omettrons cette dépendance dans ce qui suit. On pourrait fixer
$v_0=\frac{1}{p+1}$.


\vspace{2 mm}

Si $\kappa=(k,\chi) \in \WW(\CC_p)$ est de conducteur $m$, 
on dispose d'une série d'Eisenstein $E_{\kappa}$, qui est une forme modulaire de poids
$k$ sur $X_1(Np^m,d)(v)$ de caractère trivial hors de $p$, $\chi \tau^{-k}$
en $p$. Si $0\leq v<\frac{p}{p^{m-1}p+1}$, $Z_1(Np,d)(v)$ s'identifie
canoniquement au quotient de $Z_1(Np^m,d)(v)$ par l'action des diamants de $(1+p\ZZ_p)/(1+p^m\ZZ_p)$,
permettant de voir $M_0(p,d,v)$ comme le sous-espace de $M_0(p^m,d,v)$ de caractère trivial sous
$(1+p\ZZ_p)/(1+p^m\ZZ_p)$. La multiplication par $E_{\kappa}$ est un isomorphisme
de $\CC_p$-Banach de $M_0(p^m,d,v)$ sur $M_k(p^m,d,v)$ multipliant le caractère en $p$ par $\chi \tau^{-k}$. Le $q$-développement à l'infini de $E_{\kappa}$ est la
spécialisation en $\kappa$ d'un $q$-développement "abstrait" $\EEE(q) \in 1+q\Lambda[[q]]$
appelé famille d'Eisenstein restreinte (\cite{eigen}, \S 2.2). \par
\vspace{2 mm}

On pose $F:=(F(V,n))$, le système de modules de Banach sur $\WW$ associé au
système d'espaces de Banach $\{ A(Z_1(Np,d)(v_n)), \res_n, n\in \NN\}$, les
applications $$\res_n: A(Z_1(Np,d)(v_n)) \rightarrow
A(Z_1(Np,d)(v_{n+1}))$$ étant les restrictions compactes naturelles.  
Par définition, si $V \subset \WW$ est ouvert affinoide, $F(V,n)=A(V \times
Z_1(Np,d)(v_n))$. $F$ est de manière
naturelle une représentation de $\HH$, $U_p \in \Comp(F)$ (\cite{eigen}
3.2, \cite{BMF} B5).  Si $\kappa=(k,\chi) \in \WW(\CC_p)$ est de conducteur
$m$, $n\geq m-1$, $F_{\kappa,n}$ s'identifie, 
comme représentation de $\HH$, au sous-espace de $M_k(p^m,d,v_n)$ de caractère $\chi
\tau^{-k}$ sous $(\ZZ/p^m\ZZ)^*$. \par \vspace{2 mm}
En chaque pointe $c$ dans $Z_1(Np,d)(0)$, la théorie des courbes de
Tate (voir \cite{K} A1.2) nous fournit un $q_c$-développement: $F(V,n) \hookrightarrow
A(V)[[q_c]]$. 
De plus, $(\ZZ/Np\ZZ)^*$ agit par automorphismes naturels sur
$X_1(Np,d)$ en préservant $Z_1(Np,d)(v)$ et donc l'ensemble de ses pointes. 
Notons $A(Z_1(Np^m,d)(v_n))^{0,\varepsilon}$ le sous-espace de
$A(Z_1(Np^m,d)(v_n))$ composé des fonctions de $q_c$-développement nul pour tout $c$ dans
$Z_1(Np,d)(0)$, et sur lequel $(\ZZ/Np\ZZ)^*$ agit par le caractère $\varepsilon$ fixé plus haut.
On s'intéresse au système de modules de Banach $F^{0,\varepsilon}$ associé au
système d'espaces de Banach $$\{ A(Z_1(Np,d)(v_n))^{0,\varepsilon},
{\res_n}_{|A(Z_1(Np,d)(v_n))^{0,\varepsilon}}, n \in \NN\}$$ C'est un
sous-système de modules de Banach de $F$, préservé par $\HH$.
Son évaluation $F^{0,\varepsilon}_{\kappa,n}$ en $\kappa=(k,\chi) \in
\WW(\CC_p)$ de conducteur $m$ tel que $n \geq m-1$ s'identifie au sous-espace de
$M_k(p^m,d,v_n)$ composé des formes de caractère $\varepsilon \tau^{-k} \chi $ s'annulant
aux pointes de $Z_1(Np^m,d)(0)$ (noter que $E_\kappa$ est non nulle en
chacune de ces pointes).

\par \vspace{2 mm}
	On définit dans ce qui suit le sous-module de Banach de
$F^{0,\varepsilon}$ composé des familles de formes surconvergentes "nouvelles en
$d$" (voir aussi \cite{BMF} B5). 
Soit $l$ premier divisant $d$, $ld'=d$, on dispose de morphismes canoniques
finis et plats $$\pi_l:  X_1(Np^m,d) \rightarrow X_1(Np^m,d')$$ oubliant le sous-groupe d'ordre $l$.
Ces morphismes induisent des morphismes rigides analytiques 
$\pi_{l}: Z_1(Np^m,d)(v) \rightarrow Z_1(Np^m,d')(v)$, finis et plats de degré $l+1$ (étale hors des
pointes). L'isomorphisme canonique $\pi_l^*(\omega)\simeq \omega$ permet d'en
considérer la trace, ${\pi_l}_*\omega \rightarrow \omega$, d'où en
particulier sur les sections sur $Z_1(Np^m,d')(v)$:
$$tr(\pi_l)_k: M_k(p^m,d,v) \rightarrow M_k(p^m,d',v)$$ 

On étend $tr(\pi_l)_0: A(Z_1(Np,d)(v)) \rightarrow
A(Z_1(Np,d')(v))$, linéairement en un $A(V)$-morphisme 
$$\Tr_l: A(V\times Z_1(Np,d)(v)) \rightarrow A(V \times Z_1(Np,d')(v))$$ \par \vspace{2 mm}

$\Tr_l$ définit ainsi un endomorphisme du module de Banach
$F^{0,\varepsilon}$, il est bien défini pour tous les couples
$(V,n)$.

%
%
%
%
%

\begin{lemme} Soit $\kappa=(k,\chi) \in \WW(\CC_p)$, $m=\textrm{cond}(\kappa)$, 
$n\geq m-1$, le diagramme suivant est commutatif:

$$\xymatrix{ F_{\kappa,n} \ar@{->}[d]_{\Tr_l} \ar@{->}[r]^{\hspace{- 1 cm}\times E_{\kappa}}
& M_k(p^m,d,v_n) \ar@{->}[d] \ar@{->}[d]^{tr(\pi_l)_k} \\
 F_{\kappa,n} \ar@{->}[r]^{\hspace{-1 cm} \times E_{\kappa}}
& M_k(p^m,d',v_n) }$$

\end{lemme}

{\it Preuve:} Il faut montrer que si $f \in M_0(p^m,d,v)$,
$tr(\pi_l)_k(E_{\kappa}f)=E_{\kappa}tr(\pi_l)_0(f)$, ce qui découle
immédiatement de ce que $E_{\kappa}$ est de niveau premier à $l$. $\square$ 

\par \vspace{2 mm}

On vérifie aisément que $\Tr_l$ commute aux correspondances de Hecke hors de
$l$, aux opérateurs diamants, et que $\Tr_l^2=(l+1)\Tr_l$ (cf. \cite{BMF}
B5.1). Notons $W_l$ l'involution d'Atkin sur $X_1(Np^m,d)/\ZZ_p$ définie
modulairement par $(E,H,\alpha) \mapsto (E/H,E[l]/H,\pi\cdot \alpha)$, 
$H$ étant un sous-groupe d'ordre $l$ de $E$, $\pi$ l'isogénie $E \rightarrow E/H$, et $\alpha$ une structure de niveau de type $\Gamma_1(Np^m)\cap
\Gamma_0(d')$. Elle préserve les $Z_1(Np^m,d)(v)$ et induit par extension des
scalaires une involution encore notée $W_l$ sur les $A(Z_1(Np,d)(v)\times
V)$.\par \vspace{2 mm}

{\bf Définition:} Soit $\kappa \in \WW(\CC_p)$, $f \in
F_{\kappa,n}^{0,\varepsilon}$ est dite nouvelle en $d$ si pour tout $l$
divisant $d$, $\Tr_l(f)=\Tr_l(W_l(f))=0$. \par \vspace{2 mm} 

On notera $F^{0,\varepsilon,d}$ le système de modules de Banach
sur $\WW$ associé au système d'espaces de Banach $ \{
A(Z_1(Np,d)(v_n))^{0,\varepsilon,d}, res_n , n \in \NN \}$, où 
$$A(Z_1(Np,d)(v_n))^{0,\varepsilon,d}:=A(Z_1(Np,d)(v_n))^{0,\varepsilon}\cap (\bigcap_{l|d} \Ker(\Tr_l)\cap 
\Ker(\Tr_l \cdot
W_l))$$

\begin{lemme}$F^{0,\varepsilon,d}$ est un sous-système de modules de Banach 
de $F^{0,\varepsilon}$ stable par l'action de $\HH$. 
\end{lemme}
{\it Preuve}: Ça n'est pas complètement évident en ce qui concerne les $U_l$ avec
$l|d$. Il suffit de vérifier que pour $\kappa=(k,\chi) \in \WW(\CC_p)$,
$m=\cond(\kappa)$, $n \geq m-1$, $U_l$ préserve $F_{\kappa,n}^{0,\varepsilon,d}$ dans
$F_{\kappa,n}^{0,\varepsilon}$. Sur $M_k(p^m,d,v_n)$, on dispose d'un 
endomorphisme $W_{l,k}$ défini par
$W_{l,k}(f)(E,H,\alpha,\omega)=f(E/H,E[l]/H,\pi.\alpha,
(\pi^{\vee})^*(\omega) )$,
avec les notations évidentes. On vérifie que sur le sous-espace de caractère
$\varepsilon\chi\tau^{-k}$, on a $W_{l,k}^2=\varepsilon(l)\kappa(l)$ et $$tr(\pi_l)_k(f)=f+l\varepsilon(l)^{-1}\kappa(l)^{-1}U_l(W_{l,k}(f)), \, \,
tr(\pi_l)_k(W_{l,k}(f))=f+l\cdot U_l(f)$$ De plus, si $g \in
F_{\kappa,n}^{0,\varepsilon}$ et $e_l$ désigne la fonction inversible sur
$X_1(Np,d)(v_n)$ de $q$-développement $\frac{E_{\kappa}(q^l)}{E_{\kappa}(q)}$, alors
un calcul montre que
$$E_{\kappa}^{-1}W_{l,k}(gE_{\kappa})=e_lW_l(g)$$
En particulier, la multiplication par $E_{\kappa}$ envoie
$F_{\kappa,n}^{0,\varepsilon,d}$ isomorphiquement sur l'intersection des $\Ker(\tr(\pi_l)_k)\cap
\Ker(\tr(\pi_l)_k\cdot W_{l,k})$ dans le sous-espace de $M_k(p^m,d,v_n)$ de caractère
$\varepsilon\chi\tau^{-k}$. Les relations ci-dessus montrent que ce
sous-espace est stable par $W_{l,k}$ puis par $U_l$, qui coïncide avec $-l^{-1}W_{l,k}$ sur ce dernier. $\square$ \par \vspace{2 mm}

Si $\kappa=(k,\chi) \in \WW(\CC_p)$, $n\geq \cond(\chi)-1$, $F^{0,\varepsilon,d}_{\kappa,n}$ est
$\CC_p\otimes_{\ZZ}\HH$-isomorphe au sous-espace de Banach de 
$M_k(p^m,d,v_n)$ constitué des formes s'annulant aux pointes de
$Z_1(Np^m,d)(0)$, de caractère $\varepsilon \chi \tau^{-k}$, et annulées par
les $\tr(\pi_l)_k$ et $\tr(\pi_l)_k\cdot W_{l,k}$. Ces formes seront
appelées {\it nouvelles en $d$}. Sur le sous-espace de
$M_k(p^m,d,v_n)$ composé des restrictions à $Z_1(Np^m,d)(v_n)$ des formes
convergentes sur tout $X_1(Np^m,d)$, cette condition d'être nouvelle en $d$
est précisément la condition usuelle. \vspace{3 mm}

{\it Remarques}:  i) Soit $\kappa=(k,\chi) \in \WW(\CC_p)$ de conducteur $m$, $r \in \NN$, $f \in
F_{\kappa,r}^{0,\varepsilon,d}$ propre pour tout $\HH$, on note $\chi: \HH \rightarrow \CC_p$ le caractère
défini par $T_n(f):=\chi(T_n)f$. $f$ a un $q$-développement sur la pointe $\infty$ de la
forme $\sum_{n\geq 1} a_n q^n$ tel que $a_n=\chi(T_n)a_1$. Le principe du
$q$-développement (i.e l'irréductibilité de $Z_1(Np^m,d)(0)$) montre alors que $a_1 \neq 0$, et donc que l'on peut
supposer $a_1=1$. Ceci définit donc une bijection entre caractères de $\HH$
dans $F_{\kappa,r}^{0,\varepsilon,d}$ et les éléments de ce dernier qui sont propres et de
$q$-développement normalisé à $1$ ("multiplicité $1$ faible"). \par
\vspace{1 mm}
ii) Les $\CC_p$-espaces de Banach $A(Z_1(Np,d)(v))$, ainsi que tous leurs sous-espaces
considérés dans ce paragraphe, sont orthonormalisables. En effet, d'après
\cite{Ser} 1.1, tout espace de Banach provenant par extension des
scalairs d'un espace de Banach sur un corps local est orthonormalisable. De
plus, par \cite{Ser} 1.2, les inclusions $A(Z_1(Np,d)(v_n)) \supset
A(Z_1(Np,d)(v_n))^{0,\varepsilon} \supset
A(Z_1(Np,d)(v_n))^{0,\varepsilon,d}$ sont scindées dans la catégorie des
$\CC_p$-espaces de Banach. \par \vspace{2 mm}

\subsection{Formes modulaires quaternioniques $p$-adiques} ${}^{}$
\vspace{2 mm}

Dans ce qui suit, nous nous référerons à \cite{Buz2}. \par \vspace{2 mm}

\subsubsection{Séries principales $p$-adiques de l'Iwahori} \label{rep}
Introduisons tout d'abord quelques notations de théorie des groupes: 
$$L \textrm{ désigne le Borel inférieur de }  GL_2(\QQ_p)$$
$$N \textrm{  les unipotents supérieurs de }  GL_2(\ZZ_p) $$
$$I(m) \textrm{ le sous-groupe de $GL_2(\ZZ_p)$ composé des éléments triangulaires
supérieurs modulo $p^{m}$ }$$  
$$ I:=I(0) \textrm{  l'Iwahori, }\, \, \, u:=\diag(1,p)$$ 
$$\TTT(m) \textrm{ est le sous-monoide de $GL_2(\QQ_p)\cap M_2(\ZZ_p)$ engendré par
$I(m)$ et $u$}, \, \, \, \TTT:=\TTT(1)$$ 

La décomposition d'Iwahori s'écrit $I=(L\cap I)\times N$. On identifie $N$ à $\ZZ_p$ via
 $$ \left( \begin{array}{cc} 1 & t \\ 0 & 1 \end{array} \right) \mapsto t $$
La grosse cellule de Bruhat $L\backslash LI \subset
L\backslash GL_2(\QQ_p)=\PP^1(\QQ_p)$ est stable par multiplication à droite par
$\TTT$. On note $T$ la coordonnée sur $\ZZ_p$, l'action de $\TTT$ par translation à droite sur les fonctions sur $L\backslash
LI =\ZZ_p$ s'écrit alors 
$$ \left( \begin{array}{cc} a & b \\ c & d \end{array} \right).T = \frac{b+dT}{a+cT} $$
\par \vspace{2 mm}
	Soit $n \in \NN$, on note $\C_{n}$ la $\CC_p$-algèbre de Banach des fonctions sur $\ZZ_p$ de restriction analytique aux
$a+p^n\ZZ_p$, $\C_n$ est stable par l'action de $\TTT$ et $\C_0$
s'identifie à l'algèbre de Tate $\CC_p\!\!<\!T\!>$. Les restrictions naturelles $i_n: \C_n \rightarrow
\C_{n+1}$ sont compactes injectives , ce qui fait de
$\C=\{\C_n,i_n, n\in \NN\}$ un système d'espaces de Banach tel que $u \in
\Comp(\C)$. On note  $\C_{\WW}$ le système de modules de Banach sur $\WW$ associé à $\C$.
Il sera commode de poser $\C_{-n}:=\C_0$ et $\C(V,-n):=\C(V,0)$ si $n \in
\NN$, $V \subset \WW$ ouvert affinoide.  \par

\vspace{2 mm}

	Soit $j: I \rightarrow \C_0^*$ le 1-cocycle défini par  
$$j( \left( \begin{array}{cc} a & b \\ c & d \end{array} \right) ):=
\tau^{-1}(a)(a+cT) \in 1\!+\!p\ZZ_p\!<\!T\!> \, \subset \, \C_0^*$$ Il se prolonge à
$\TTT$ en le prenant trivial sur $u$.  Si $\kappa \in \WW(\CC_p)$, $\gamma \in
\TTT$, on définit un cocycle $\kappa(j): \TTT \rightarrow \bigcup_{n \in \NN}
\C_{n}^*$ par
$$(\kappa(j)(\gamma))(t):=\kappa(j(\gamma)(t)), \, \, \, t \in \ZZ_p$$
Si $\gamma \in \TTT(m)$ et $\kappa \in \WW_r(\CC_p)$, $\kappa(j)(\gamma) \in
\C_{r-m}^*$.
On notera $\rho_{\kappa}$ la représentation de $\TTT$ sur
l'espace $\bigcup_{n \in \NN}\C_n$ tordue par $\kappa(j)$, i.e
$\rho_{\kappa}(v):=\kappa(j)(\gamma)\gamma.v$. Cette torsion n'affecte pas
$u$ et si $\kappa \in \WW_r(\CC_p)$, $\TTT(m)$ préserve $\C_n$ dès que $n\geq
r-m$. Plus généralement, si $V \subset \WW$ est un ouvert
affinoide, on dispose d'un $1$-cocycle $\kappa^{univ}(j): \TTT \rightarrow
\bigcup_{n \in \NN}\C(V,n)^*$, tel que 
$\kappa^{univ}(j)_{\kappa}=\kappa(j)$ si $\kappa \in \WW(\CC_p)$ et 
$\kappa^{univ}(j)(\gamma) \in \C_{r-m}^*$ si $\gamma \in \TTT(m)$ et 
$V \subset \WW_r$. On dispose donc d'une représentation
$\rho^{univ}$ de $\TTT$ sur $\bigcup_{n \in \NN} \C(V,n)$, en
tordant par $\kappa^{univ}(j)$ la représentation naturelle obtenue par
extension complète des scalairs à $A(V)$ de celle sur $\C_n$. Si $V \subset
\WW_r$, $\TTT(m)$ préserve $\C(V,n)$ dès que $n \geq r-m$. \par \vspace{2 mm}
	Ainsi, le système de modules de Banach $\C_{\WW}$ sur
$\WW$ est muni d'une représentation $\rho^{univ}: \TTT \rightarrow
\End(\C_{\WW})$. Si $V \subset \WW_r$, $\TTT(m)$ préserve $\C(V,n)$ dès que
$n \geq r-m$, et on a $u \in \Comp(\C_{\WW})$. On appelle $\C_{\WW}$ {\it la famille analytique des séries principales
$p$-adiques de $I$}. \par \vspace{2 mm}

{\it Remarques:}  i) $\C^{\dagger}=\bigcup_{n \in \NN}\C_n$ est
l'espace des fonctions localement analytiques sur $\ZZ_p$, muni de sa
topologie localement convexe c'est un espace de type compact au sens de
\cite{ST} \S 1 \footnote{Strictement, il faudrait plutôt prendre des fonctions
$K$-valuées avec $K$ sphériquement complet.}. La représentation $\rho_{\kappa}$ de $I$ sur $\C^{\dagger}$ est la série principale localement
analytique de $I$ de caractère $\kappa$ (\cite{ST} \S 5, noter que leur $B$ est
l'Iwahori opposé à $I$).\par 
ii) Les assertions de ce paragraphe sont détaillées dans \cite{Buz2} \S
4, \S
7, noter qu'il a des actions à droites, non à gauche. A
cette modification près, on a dans ses notations:
$\mathcal{A}_{\kappa,p^{-n}}=\C_{\kappa,n}$ et $M_{m} \supset \TTT(m)$; d'autre part
si $\kappa \in \WW_r(\CC_p)\backslash \WW_{r-1}(\CC_p)$, "$m$ is
good for $(\kappa,p^{-n})$" équivaut à $n \geq r-m$. 

\subsubsection{Formes modulaires} \label{quat}

Soit $D(\QQ)$ une algèbre de quaternions sur $\QQ$, on fixe
$D(\ZZ)$ un ordre maximal de $D(\QQ)$, et on note
$D$ le schéma en anneaux sur $\ZZ$ associé, $D^*$ son groupe des inversibles.
On suppose que $D$ est définie. Soit $d=\disc(D)$ le produit des premiers ramifiés,
on fixe un isomorphisme au dessus de $\ZZ[1/d]$, 
$$\varphi: D/\ZZ[1/d] \simeq \mathbb{M}_2/\ZZ[1/d]$$ Si $M$ est un entier
premier à $d$, on notera $U_1(M)$ (resp. $U_0(M)$) le sous-groupe
ouvert compact de $D^*(\AAA_f) \simeq D^*(\AAA_{f}^{(d)}) \times \prod_{l | d}
D^*(\QQ_l)$, décomposé selon ce produit, valant le compact maximal $D^*(\ZZ_l)$ en les premiers $l$
divisant $d$, égal à $\varphi^{-1}(\Gamma_1(M))$ (resp.
$\varphi^{-1}(\Gamma_0(M))$) sur l'autre facteur, avec:

$$ \Gamma_1(M)=\{ g \in GL_2(\widehat{\ZZ[1/d]}), \, \, g \equiv \left(
\begin{array}{cc} 1 & * \\ 0 & * \end{array} \right)\,\bmod M\},$$ 
$$\Gamma_0(M)=\{ g \in GL_2(\widehat{\ZZ}), \, \, g \equiv
\left( \begin{array}{cc} *
& * \\ 0 & * \end{array} \right)\,\bmod M \}$$

On verra les caractères de $(\ZZ/M\ZZ)^*$ comme des caractères de $U_0(M)$
par l'étoile supérieure gauche. Si $M\geq 5$, $D^*(\QQ) \times
U_1(M)$ agit librement sur $D^*(\AAA_f)$, et $D^*(\QQ)\backslash
D^*(\AAA_f)/U_1(M)$ est fini de cardinal noté $h_1(M)$. On choisit $M=Np$ comme en \ref{formes}. \par \vspace{1 mm}
On note $\FF^{D}$ le système de modules de Banach sur $\WW$ tel que si $V
\subset \WW$ est ouvert affinoide, $n \in \NN$, 
$F^{D}(V,n)$ est le $A(V)$-module de Banach des fonctions 
$f: D^*(\QQ)\backslash D^*(\AAA_f)\rightarrow \C(V,n)$ satisfaisant 
$$f(xu)=j(u_p^{-1})^{-2}\rho^{univ}(u_p^{-1})f(x), \, \,  \forall (x,u) \in D^*(\AAA_f) \times
U_1(Np) $$ et $i_n$ est l'application déduite de la restriction canonique
$\C(V,n)
\rightarrow \C(V,n+1)$. L'orthonormalisabilité vient de ce que $F^D(V,n)$ 
s'identifie à $\C(V,n)^{h_1(Np)}$ car $Np\geq 5$ (voir aussi \cite{Buz2} \S 7, \S 4). \par
\vspace{2 mm}
	On a une représentation naturelle sur $F^D$ de $U_0(Np)/U_1(Np)\simeq
 (\ZZ/Np\ZZ)^*$, définie par $(<\gamma>.f)(x)=j(\gamma)^{-2}\rho^{univ}(\gamma)f(x\gamma)$.
On notera $F^{D,\varepsilon}$ le sous-système de modules de Banach de $F^D$ sur
lequel $(\ZZ/Np\ZZ)^*$ agit par $\varepsilon^{-1}$. On dispose d'une
représentation naturelle $\HH \rightarrow \End(F^{D,\varepsilon})$ telle que
$U_p \in \Comp(F^{D,\varepsilon})$. Les doubles classes considérées ici pour
les opérateurs de Hecke sont (avec les abus évidents) les $\diag(1,l)$, pour $T_l$ et $U_l$
si  $l\nmid d$, $\diag(l,l)$ pour $S_l$ si $l \nmid Npd$, comme en \ref{JLclass} pour $U_l$
avec $l|d$. Si $\kappa \in \WW(\CC_p)$, $F^{D,\varepsilon}_{\kappa,n}$
s'identifie à l'espace des fonctions $D^*(\QQ)\backslash D^*(\AAA_f)
\rightarrow \C_n$ telles que 
$$f(xu)=\varepsilon^{-1}(u)j(u_p^{-1})^{-2}\rho_{\kappa}(u_p^{-1})f(x), \,
\,\, \forall \, (x,u) \in D^*(\AAA_f) \times U_0(Np)$$

$(F^{D,\varepsilon}_{\kappa})^{\dagger}$, vu avec sa structure de
$\HH$-module, est {\it l'espace des formes modulaires $p$-adiques quaternioniques de poids-caractère
$\kappa$, de niveau modéré $N$, de caractère $\varepsilon$}.
\par \vspace{2 mm}

Soient $n,\, m,\, r$ des entiers tels que $n\geq r-1 \geq 0$, $n\geq m-1 \geq
0$, et $\kappa \in \WW_r(\CC_p)$, considérons l'espace
annexe $F_{\kappa,n}^{D,\varepsilon}[m]$ des fonctions $f:
D^*(\QQ)\backslash D^*(\AAA_f)\rightarrow
\C_{\kappa,n}$ satisfaisant
$$f(xu)=\varepsilon^{-1}(u)j(u_p^{-1})^{-2}\rho_{\kappa}(u_p^{-1})f(x), \, \,
\forall (x,u) \in
D^*(\AAA_f^{(d)}) \times U_0(Np^m) $$

C'est un $\HH$-module, en prenant cette fois-ci des doubles classes par rapport à
$U_0(Np^m)$ (ce qui ne change que $U_p$, encore défini par la double classe de
$\diag(1,p)$). D'après \cite{Buz2} (\S 5 lemme 3,iv), 
$F_{\kappa,n-m+1}^{D,\varepsilon}[m]$ et
$F_{\kappa,n}^{D,\varepsilon}$ sont isomorphes comme $\HH \otimes_{\ZZ}\CC_p$-modules. 
Si de plus $\kappa=(k,\chi)$ est de conducteur $m$, et $k\geq 2$,
$F_{\kappa,0}^{D,\varepsilon}[m]$ contient comme sous-$\HH\otimes_{\ZZ}\CC_p$-module 
l'espace des fonctions à valeurs polynomiales en $T$ de degré $\leq k-2$, ce dernier est 
isomorphe comme $\HH \otimes_{\ZZ} \CC_p$-module à l'espace des fonctions 
$f: D^*(\QQ)\backslash D^*(\AAA_f)\rightarrow \textrm{Sym}^{k-2}(\CC_p^2) $ telles
que $$  f(xu)=\varepsilon(u)^{-1}\chi(u_p)^{-1}\tau^{k}(u_p)u_p^{-1}f(x), \, \,  \forall
\, \, (x,u) \, \, \in \, \, D^*(\AAA_f) \times
U_0(Np^m)$$

Ce $\HH \otimes_{\ZZ} \CC_p$ module est $S_{k}^D(Np^m,\varepsilon\chi\tau^{-k},\CC_p)$
dans la notation du \S \ref{pstructure}.

\vspace{2 mm}

\section{Préliminaires de théorie spectrale}


\subsection{Semi-simplification en dimension infinie} \label{semi} Soit $K$ un corps
valué complet non archimédien (non discret), $V$ un $K$-espace de Banach orthonormalisable, $U$ un endomorphisme compact de $V$. 
La série caractéristique $P(T)=\det(1-TU)$ de $U$ se décompose sous la forme $$P(T)=\prod_{i \geq 0}
P_i(T)^{n_i}$$ où les
$P_i(T)$ sont des irréductibles de $1+TK[T]$ deux à deux distincts tels que $|P_i(T)-1|
\underset{i \rightarrow \infty}{\rightarrow} 0$ pour la norme $|.|$ sur
$K[T]$ du sup des coefficients. Par \cite{Ser}, on sait que $V$ est somme directe topologique de $Ker(U)$ et
des espaces de dimension finie $V(P_i):=\Ker(P^*_i(u)^{n_i})$, $Q^*(T)$
désignant le polynôme réciproque de $Q(T)$, sur lesquels $U$ a pour polynôme
caractéristique $P^*_i(T)^{n_i}$. Soit $H$ une
$K$-algèbre, $\rho: H \rightarrow \End_K(V)$ une représentation telle que
$\rho(H)$ contienne $U$ et lui commute, alors $\rho(H)$ stabilise
les $V(P_i)$, que l'on peut semi-simplifier. \par \vspace{2 mm}

{\bf Définition:}  On notera $\mathcal{X}_U(V)$ l'ensemble des
représentations irréductibles de $H$ apparaissant dans
la réunion des semi-simplifications des $V(P_i)$, comptées avec multiplicité
(qui sont nécessairement finies). 
\par \vspace{2 mm}
Notons que $\mathcal{X}_U(V)$ dépend de
$U$, mais que $\mathcal{X}_U(V)=\mathcal{X}_{U'}(V)$ si $U' \in \rho(H)$ est un
autre endomorphisme compact de $V$ commutant à $\rho(H)$ et à $U$ tel que
$\Ker(U')=\Ker(U)$. Si $\mathcal{X}$ est un ensemble de représentations d'une
algèbre, on notera $|\mathcal{X}|$ l'ensemble des classes d'isomorphie de
représentations apparaissant dans $\mathcal{X}$ (autrement dit, on oublie les
multiplicités).
\begin{prop} \label{semisimpl} Soient $(\rho_1,V_1)$ et $(\rho_2,V_2)$, des
représentations d'une $K$-algèbre $H$ dans les endomorphismes continus de $V_1$
et $V_2$, telles que $V_i$ ($i=1,2$) est muni d'un 
endomorphisme compact $U_i$ commutant à $\rho_i(H)$. \par 
On suppose de plus que pour
tout $h \in H$, $\det(1-T\rho_1(h)U_1)=\det(1-T\rho_2(h)U_2) \in K\{\{T\}\}$.
Alors $\mathcal{X}_{U_1}(V_1)=\mathcal{X}_{U_2}(V_2)$. \end{prop}
{\it Preuve:} 
Soit $\alpha \in \RR$, $i \in \{1,2\}$, $V_i^{\alpha}$ le plus grand
sous-espace de dimension finie de $V_i$ stable par $U_i$ sur lequel le
polygone de Newton du polynôme caractéristique de $U_i$ est de pente $\alpha$. 
Soit $h \in H$, $\rho_i(h)$ stabilise les $V_i^{\alpha}$, et ses valeurs
propres sur ces derniers sont toutes bornées par $||\rho_i(h)||< \infty$,
$||.||$ désignant la norme d'opérateur sur $V_i$. Soit $x \in K^*$ tel que
$|x|<1$, on peut trouver par conséquent un $N \in \NN$ tel que
$\rho_i(1+x^Nh)$ ait toutes ses valeurs propres de norme $1$ sur les
$V_i^{\alpha}$ ($i=1,2$). On pose $h'=1+x^Nh \in H$. En co-trigonalisant $U_i$ et
$\rho_i(h')$ sur $V_i^{\alpha}$, on en déduit:
$$ \det(1-TU_i\rho_i(h'))^{\alpha}=\det(1-TU_i\rho_i(h')_{|V_i^{\alpha}})$$
Si $Q \in 1+TK\{\{T\}\}$, $Q^{\alpha}$ désigne le polynôme $\in 1+TK[T]$
divisant $Q$ tel que le polygone de Newton de $Q/Q^{\alpha}$ n'a pas la
pente $\alpha$. Ainsi, $\det(1-T\rho_1(h')U_1)=\det(1-T\rho_2(h')U_2)$ donne 
$$\det(1-T\rho_1(h'){U_1}_{|V_1^{\alpha}})= \det(1-T\rho_2(h'){U_2}_{|V_2^{\alpha}})$$ 
$\rho_i(h)U_i$ étant injectif sur $V_i^{\alpha}$, notons que cela implique
que $\dim(V_1^{\alpha})=\dim(V_2^{\alpha})$, et que $\rho_1(h'){U_1}_{|V_1^{\alpha}}$
et $\rho_2(h'){U_2}_{|V_2^{\alpha}}$ ont même polynôme caractéristique. Ceci
reste vrai pour les même raisons en remplaçant $h'$ par $h'+\lambda$ pour
(une infinité de) $\lambda \in K$ assez petit. Si $A$ et $B$ sont deux
endomorphismes qui commutent d'un $K$-espace vectoriel de dimension finie
$r$, avec $A$ inversible, la donnée de
$$\det(X.1-A(B+Y.1))=\prod_{i=1}^r(X-a_iY-b_i) \in \overline{K}[X,Y]$$ 
permet de retrouver $\det(X.1-B)=\prod_{i=1}^r(X-(b_i/a_i))$ (ici
$\overline{K}$ est une clôture algébrique de $K$).
On en déduit $\det(T-\rho_1(h')_{|V_1^{\alpha}})=\det(T-\rho_2(h')_{|V_2^{\alpha}})$, puis
la même chose en remplaçant $h'$ par $h$. Ainsi, $V_1^{\alpha}$ et $V_2^{\alpha}$
sont deux représentations de $H$ ayant même polynômes caractéristiques, on
sait alors que leurs semi-simplifications sont isomorphes, ce qui conclut.
$\square$ \par \vspace{2 mm}

	Terminons cette partie par une légère amélioration de
\ref{semisimpl}. Soient $E^i$, $i=1,2$, deux systèmes d'espaces de
Banach, munis de représentations $\rho_i: H \rightarrow \textrm{End}(E^i)$ et 
d'endomorphismes compacts $U_i \in \textrm{Comp}(E^i)$, $U_i \in
\rho_i(H)$ commutant à $\rho_i(H)$. $E^{i,\dagger}$ est muni d'une
opération de $\rho_i(H)$ et $U_i$. Si $\alpha \in \RR$ est fixé et $E_n^{i,\alpha}(U_i)$ désigne le sous-espace de
$E^i_n$ sur lequel $U_i$ est de pente $\alpha$, $i_n$ induit pour tout $n$
assez grand une bijection $E_n^{i,\alpha}(U_i) \rightarrow E_{n+1}^{i,\alpha}(U_i)$ qui commute à
$U_i$. Cet espace définit donc un sous-espace de dimension finie 
$E^{i,\alpha}(U_i)$ de $E^{i,\dagger}$ qui hérite d'une représentation de $H$. Ceci permet de 
définir à nouveau $\mathcal{X}_{U_i}(E^i)$ comme étant l'ensemble des
représentations (comptées avec multiplicités) de $H$ apparaissant dans les
semi-simplifications des $E^{i,\alpha}(U_i)$, $\alpha$ variant dans $\RR$.  

\begin{cor} \label{corsemisimpl} Sous ces hypothèses, supposons que pour
tout $h \in H$, $\det(1-T\rho_1(h)U_1)=\det(1-T\rho_2(h)U_2) \in K\{\{T\}\}$.
Alors $\mathcal{X}_{U_1}(E^1)=\mathcal{X}_{U_2}(E^2)$. 
\end{cor}
{\it Preuve:} Fixons $\alpha \in \RR$, $E^{i,\alpha}(U_i)$ est de dimension
finie, $H$ agit donc sur $E^{1,\alpha}(U_1)\oplus E^{2,\alpha}(U_2)$ à
travers un quotient de dimension finie. Considérons une
sous-$K$-algèbre $H_0$ de $H$ de type fini sur $K$ engendrant
toute l'image de $H$ dans ce quotient; pour $n$ assez grand, $H_0$ et $U_i$
agissent alors par endomorphismes sur $E^i_n$, $i=1,2$. On applique
la proposition \ref{semisimpl} avec $V_i=E^i_n$, il vient
$\mathcal{X}_{U_1}(E^{1,\alpha}(U_1))=\mathcal{X}_{U_2}(E^{2,\alpha}(U_2))$,
ce qui conclut. $\square$

\subsection{Un critère de densité} \label{critere}
	On fixe $M$ un système de modules de Banach comme en
\S \ref{banachsystem}, $\textrm{dim}(\WW)>0$, muni d'une action d'une
$\CC_p$-algèbre $H$, et d'un endomorphisme compact $U$ lui commutant, on
notera $(M,H,U)$ une telle donnée. On rappelle que si $T/\CC_p$ est un espace
rigide, un sous-ensemble $X\subset
T(\CC_p)$ est dit Zariski-dense si pour tout fermé analytique $F$ de $T$ 
tel que $X \subset F(\CC_p)$, alors $F(\CC_p)=T(\CC_p)$. 
Soit $X \subset \WW(\CC_p)$ un sous-ensemble Zariski-dense, tel
que pour tout $x \in X$ et tout voisinage ouvert affinoide irréductible $V$ de $x$ dans $\WW$,
$V(\CC_p) \cap X$ est Zariski-dense dans $V$. On dira alors que $X$ est {\bf très
Zariski-dense dans $\WW$}. C'est par exemple le cas des points de la
forme $\zeta(1+p)^k$, avec $\zeta^{p^m}=1$ et $k,m \in \NN$, dans la boule
ouverte de centre $1$ de rayon $1$ de $\CC_p$. \par \vspace{1 mm}
	On se fixe de plus une "structure classique sur $X$", on entendra
par là la donnée pour tout $x \in X$ d'un sous-espace vectoriel de dimension finie $M_x^{cl}$ de
$M_x^{\dagger}$ stable par l'action de $H$. Soit $\alpha \in \RR$, on notera
$M_{x}^{\leq \alpha}$ (resp. $M_x^{\alpha}$)
le sous-espace de dimension finie de $M_x^{\dagger}$ sur lequel
$U$ est de pente au plus $\alpha$ (resp. exactement $\alpha$). On fera de plus l'hypothèse de "contrôle": 
\vspace{2 mm}
\begin{center} \label{Cl} (Cl) \hspace{3 mm} \textrm{Si $\alpha \in \RR$,
$U \subset \WW$ ouvert affinoide, 
alors pour tout
$x \in X \cap U(\CC_p)$ sauf peut-être un nombre fini d'entre eux, $M_{x}^{\leq \alpha} \subset
M_{x}^{cl}$ }
\end{center}

\vspace{1 mm}
\begin{prop} \label{det} Soient $(M_1,H,U)$ et $(M_2,H,U)$ deux systèmes de modules de Banach sur $\WW$
relativement factoriel. On se donne $X \subset \WW(\CC_p)$ un ensemble très Zariski-dense,
et une structure très classique sur $X$ pour $M_1$ et $M_2$, 
chacune de ces structures satisfaisant (Cl). Soit $h \in H$, supposons 
$$\forall x \in X, \, \, \, \, \, \det(1-ThU_{|M_{1,x}^{cl}})=\det(1-ThU_{|M_{2,x}^{cl}})\in \CC_p[T]$$ 
Alors \Fred$_{M_1}(hU)$=\Fred$_{M_2}(hU)$.
\end{prop}
{\it Preuve:} Soit $Z_i \subset \WW \times \AAA^{1}$ l'hypersurface de Fredholm de
$P_i$:=\Fred$_{M_i}(hU)$, $p_i$ la première projection $Z_i \rightarrow
\WW$. On dira que $z \in Z_i(\CC_p)$ est classique si $z=(x,\lambda)$ avec
$x \in X$ et si $\lambda^{-1}$ est valeur propre de $hU$ sur $M_{i,x}^{cl}$.
On montrera plus bas que les points classiques sont Zariski-denses dans
$Z_i(\CC_p)$, admettons le pour l'instant. 
Par hypothèse, $P_1$ s'annule sur les points classiques de
$Z_2(\CC_p)$, donc sur $Z_2(\CC_p)$ par Zariski-densité. Par symétrie,
il vient $Z_1^{\rr}=Z_2^{\rr}$ et on déduit de \cite{Con} 4.3.2 que $P_1$ et $P_2$ ont même facteurs irréductibles, il reste à prouver
que ces derniers ont même multiplicités. \par
Soit $\Pi$ un de ces facteurs irréductibles, de multiplicité $n_i$ dans $P_i$, $Z(\Pi)_i \subset Z_i$ la composante
irréductible associée. Soit $W_i$ l'ouvert de $Z(\Pi)_i$ dont le
complémentaire est l'ensemble points de $Z(\Pi)_i$ qui sont dans au moins deux
composantes irréductibles de $Z_i$. Admettons pour l'instant que l'on puisse trouver $x \in X$ et $z=(x,\lambda) \in
W_1(\CC_p)=W_2(\CC_p)$ tels que $\lambda$ soit une racine de $\det(1-ThU_{|M_{i,x}^{cl}})$ 
mais pas de $P_i(x,T)/\det(1-ThU_{|M_{i,x}^{cl}})$. 
Par le choix de $W_i$, $\lambda$ est une racine de $P_i(x,T)$ qui est en fait une racine de
$\Pi(x,T)$ mais pas des autres facteurs irréductibles. La multiplicité de
$\lambda$ comme racine de $P_i(x,T)$ est donc de la forme $n_in$ où $n$ est la multiplicité de
$\lambda$ comme racine de $\Pi(x,T)$. Mais, par le choix de $z$, $nn_i$ est aussi la
multiplicité de $\lambda$ comme racine de $\det(1-ThU_{|M_{i,x}^{cl}})$,
qui ne dépend pas de $i$. Ainsi, $n_1n=n_2n$, puis $n_1=n_2$.\par
	Il reste à trouver un tel $z$ et à prouver que les points classiques sont Zariski-denses dans $Z_i(\CC_p)$. Par hypothèse sur $\WW$
et \cite{Con} 4.3.2, les composantes irréductibles de $Z_i$ sont des
hypersurfaces de Fredholm. Ces dernières étant d'image Zariski-ouverte dans $\WW$, elles contiennent toutes des points d'image
(par $p_i$) dans $X$. Soit $z_i$ un de ces points, appartenant a une composante irréductible
$T_i$ de $Z_i$, et soit $\Omega_i \in \mathcal{C}(Z_i)$ contenant $z_i$.
$\mathcal{C}(Z_i)$ désigne le recouvrement canonique de l'hypersurface
de Fredholm $Z_i$ (voir la discussion au début de \ref{unicite}).
$\Omega_i$ étant fini et plat sur son image $V_i \subset \WW$ (que l'on peut
supposer irréductible), chacune de ses composantes irréductibles se surjecte sur $V_i$. 
$hU$ et $U$ sont des endomorphismes de $M_i(V_i,n)$ pour un $n$
assez grand que l'on fixe, et $P_i(T)_{|V_i}=\det(1-T{hU}_{|M_i(V_i,n)})\in
1+A(V_i)T\{\{T\}\}$. Par choix de $\Omega_i
\in \mathcal{C}(Z_i)$ et un théorème de Coleman (\cite{eigen} \S 7.1, \cite{BMF} A.4.3), 
on peut trouver un sous-$A(V_i)$-module $N_i$ de $M_i(V_i,n)$ 
localement libre de rang fini, stable par $U$ et $hU$, tel que 
$\Omega_i$ soit le fermé des zéros de $\det(1-{ThU}_{|N_i})$ dans $V_i \times
\AAA^1$. $hU$ est un endomorphisme inversible de $N_i$, ainsi donc que $U$, et ils sont
automatiquement continus (\cite{BGR} 3.7.3 proposition 2). En particulier, les valeurs propres de
$U^{-1}$ sur les $N_{i,x}$, $x \in V_i(\CC_p)$, sont bornées par une constante ne dépendant que
des $V_i$. Ceci et $(Cl)$ impliquent que pour tout $x \in X \cap V_i(\CC_p)$ sauf
peut-être un nombre fini d'entre eux, $\Omega_i(\CC_p) \cap p_i^{-1}(\{x\})$ est composé de points classiques. 
$X\cap V_i(\CC_p)$ étant Zariski-dense dans $V_i(\CC_p)$, et $\Omega_i
\rightarrow V_i$ étant fini et plat, on en déduit que les points classiques 
sont Zariski-denses dans chaque composante irréductible de $\Omega_i$, et
en particulier dans $Z_i(\CC_p)$ et $T_i(\CC_p)$ (\cite{Con} 2.2.3). Plus
exactement, on en déduit que l'ensemble des points $z=(x,\lambda)$ tels que 
$M_{i,x}^{hU=\lambda^{-1}} \subset M_{i,x}^{cl}$ est Zariski-dense dans
$T_i(\CC_p)$. $M_{i,x}^{hU=\lambda^{-1}}$ désigne ici le sous-espace de
$M_{i,x}$ qui est l'espace caractéristique pour la valeur propre $\lambda^{-1}$ 
de l'endomorphisme compact $hU$. \par 
	On appliquant cela à $T_i:=Z(\Pi)_i$. $\Omega_i \cap W_i$ est un
ouvert de $\Omega_i$ et contient donc des points du type précédent. Si
$(x,\lambda)$ en est un, $\lambda$ est une racine de $\det(1-ThU_{|M_{i,x}^{cl}})$ 
mais pas de $P_i(x,T)/\det(1-ThU_{|M_{i,x}^{cl}})$, c'est juste ce qui nous
manquait pour conclure. $\square$

\section{La correspondance de Jacquet Langlands "$p$-adique"}

\subsection{Rappels sur la correspondance classique}\label{JLclass} 	Si $M \in \NN$, on pose :

	$$ K_1(M)=\{ g \in GL_2(\widehat{\ZZ}), \, \, g \equiv \left( \begin{array}{cc} * & * \\ 0 & 1 \end{array}
\right)\,\bmod M\},$$ $$K_0(M)=\{ g \in GL_2(\widehat{\ZZ}), \, \, g \equiv \left( \begin{array}{cc} *
& * \\ 0 & * \end{array} \right)\,\bmod M \}$$

On pourra voir, par l'étoile inférieure droite, les caractères complexes de
$(\ZZ/M\ZZ)^*$ comme des caractères de $K_0(M)$. On fixe $\varepsilon$ un
tel caractère, ainsi qu'un entier $k \geq 2$. \par \vspace{1 mm}	
Le $\CC$-espace vectoriel $S_k(M,\varepsilon)$ des formes modulaires paraboliques de
poids $k$, de niveau $M$, de caractère $\varepsilon: (\ZZ/M\ZZ)^*
\rightarrow \CC$ s'identifie à l'espace des fonctions complexes sur $GL_2(\QQ)\backslash GL_2(\AAA)/K_1(M)$,
de caractère central $|.|^{2-k}\varepsilon^{-1}$, de caractère $\varepsilon^{-1}$
sous $K_0(M)$, satisfaisant les conditions
usuelles de croissance aux pointes et holomorphie, en associant à $f$ (voir par exemple \cite{hida}
\S 3.1 pour des détails) \footnote{$f_{|k g}(z):=f(g(z))j(g,z)^{-k}\det(g)^{k/2}, \, \, g \in GL_2(\RR)^+$, en
particulier, si $g \in Z(\RR)$, $f_{|k g}=f$}
$$g=(g_{\infty},g_f) \in GL_2^{+}(\RR) \times K_0(M) \mapsto
|\det(g)|^{1-k/2}f_{|_k g_{\infty}}(i)\varepsilon^{-1}(g_f)$$

Cette identification préserve l'action des opérateurs de Hecke usuels (non
renormalisés du côté adélique). Soit $l$ premier, l'opérateur de Hecke $T_l$ si $(l,M)=1$, $U_l$ sinon, est
donné par la double classe $U_0(M)\diag(l,1)U_0(M)$; si
$(l,M)=1$, $S_l=l^{k-2}\varepsilon(l)$ est l'action de $\diag(l,l)$ (ici, on voit
$\diag(a,b) \in GL_2(\AAA)$ partout trivial sauf en $l$ où il vaut
effectivement $\diag(a,b)$). \par \vspace{1 mm} 
	On notera $S_k(M,\varepsilon)^{d-new}$ le sous-espace de
$S_k(M,\varepsilon)$ composé des formes $d$-nouvelles au sens usuel. \par
\vspace{1 mm}

On fixe un entier $d$ sans facteur carré, ayant un nombre
impair de diviseurs premiers, tel que $d\divi M$ et que $\varepsilon$
soit trivial sur $(\ZZ/d\ZZ)^*$. On reconsidère $D(\QQ)$ l'algèbre de quaternions de discriminant
$d$ introduite en \ref{quat}. $K_1(M)$ (resp. $K_0(M)$) est le compact ouvert de
$D^*(\AAA_f)$ décomposé place par place, qui vaut $D(\ZZ_p)^*$ aux places
$p|d$, $\varphi^{-1}(K_1(M))$ (resp. $\varphi^{-1}(K_0(M))$) hors de $d$ (cf.
\ref{quat}), noter que $K_1(M) \neq U_1(M)$. \par \vspace{1 mm}

On note $S_k^{*,D}(M,\varepsilon)$ le $\CC$-espace vectoriel des
fonctions complexes sur $D^*(\QQ)\backslash D^*(\AAA)/K_1(M)$, de caractère
$\varepsilon^{-1}$ sous $K_0(M)$, de caractère
central $|\det(g)|^{2-k} \varepsilon^{-1}$, engendrant sous $D^*(\RR)$ un
multiple du dual de la représentation $\textrm{Sym}^{k-2}(\CC^2)$. Si $k=2$ et
$\varepsilon=1$, on note
$S$ la droite des fonctions constantes sur $D^*(\AAA)$ dans
$S_2^{*,D}(M,1)$, elle est stable par $D^*(\AAA)$. \par \vspace{2 mm}

\begin{thm}(Arthur, Jacquet, Langlands) \label{JL} Si $k>2$ ou $\varepsilon\neq 1$, les $\CC$-espaces vectoriels
$S_k^{*,D}(M,\varepsilon)$ et $S_k(M,\varepsilon)^{d-new}$ sont 
isomorphes en tant que modules sous l'algèbre de Hecke de
$GL_2(\AAA_f^{(d)})$. Si $k=2$, $\varepsilon=1$, c'est encore vrai en
remplaçant $S_2^{*,D}(M,1)$ par $S_2^{*,D}(M,1)/S$.
\end{thm}
	
	L'action des opérateurs de Hecke dans cette correspondance se précise de plus en $l|d$,
en faisant correspondre à $U_l$, l'opérateur de Hecke de $S_k^{*,D}(M,\varepsilon)$ donné par la double classe $K_0(M)\pi_lK_0(M)$,
où $\pi_l \in D^*(\AAA_f)$ est partout trivial sauf en $l$ où il vaut une
(quelconque) uniformisante de $D(\ZZ_l)$. On notera encore $U_l$ cet
opérateur de Hecke pour $D^*$. \par

	Il se trouve que nous n'allons pas considérer exactement l'espace
$S_k^{*,D}(M,\varepsilon)$, mais un autre légèrement différent (ce qui
explique la notation $*$ provisoire), qui lui est isomorphe comme module
sous l'algèbre de Hecke. Soit $$\omega_M:= \left(
\begin{array}{cc} 0 & 1 \\ M & 0 \end{array} \right) \in GL_2(\QQ)$$ on le
voit comme un élément de $D^*(\AAA)$ trivial aux places divisant $d$ et à
l'infini, diagonalement $\omega_M$ dans
$D^*(\AAA_f^{(d)})=GL_2(\AAA_f^{(d)})$. $\omega_M$ normalise
$K_0(M)$ et agit par $\diag(a,b) \mapsto \diag(b,a)$ sur le tore
diagonal de $GL_2(\AAA_f^{(d)})$. \par \vspace{2 mm}

L'application $f \mapsto \omega_M.f $ induit un isomorphisme $\CC$-linéaire de
$S_k^{*,D}(M,\varepsilon)$ sur l'espace des fonctions complexes sur
$D^*(\QQ)\backslash D^*(\AAA)/U_1(M)$
de caractère $\varepsilon^{-1}$ sous $U_0(M)$, de caractère central
$|\det(g)|^{2-k}\varepsilon^{-1}$, engendrant sous $D^*(\RR)$ un multiple du dual de
$\textrm{Sym}^{k-2}(\CC^2)$. \par \vspace{2 mm} 
On note $S^{D}_k(M,\varepsilon)$ ce $\CC$-espace
vectoriel muni de l'action de l'algèbre de Hecke obtenue par transport.
Explicitement, si $l$ est premier
ne divisant pas $d$, l'opérateur de Hecke $T_l$ si
$(l,M)=1$, $U_l$ sinon, est donné par la double classe $U_0(M)\diag(1,l)U_0(M)$; si
$(l,M)=1$, $S_l=l^{k-2}\varepsilon(l)$ est l'action de $\diag(l,l)$, $U_l$ est
inchangé si $l|d$. \par
\vspace{1 mm}

	Comme en \ref{formes}, $\HH$ désigne la $\ZZ$-algèbre de polynômes sur les lettres
$T_l$, $S_l$ si $l$ premier ne divise pas $M$, $U_l$ si $l|M$. Par ce que
l'on vient de dire plus haut, $S_k^D(M,\varepsilon)$ et
$S_k(M,\varepsilon)^{d-new}$ sont deux $\HH$-modules isomorphes, avec la
même exception pour $k=2$ que dans \ref{JL}. 

\subsection{Structures $p$-adiques des espaces de formes classiques}
\label{pstructure}
${}^{}$\par \vspace{1 mm}

Nous allons introduire une $\QQ_p$-structure sur les espaces $S_k^{D}(M,\varepsilon)$ et $S_k(M,\varepsilon)$. On fixe
pour cela $p$ premier, $\iota: \CC_p \rightarrow \CC$ un
isomorphisme de corps, $K \subset \CC_p$ un sous-corps complet, et $M=Npd$ avec $N,p$ et $d$ comme en \S   \ref{formes}.
\par \vspace{1 mm}
La courbe modulaire $X_1(M)$ a une structure naturelle sur $\QQ_p$,
préservée par les correspondances de Hecke. Soit $S_k(M,\varepsilon,K)$
le sous $K$-espace vectoriel de $H^0(X_1(M)/K,\omega^k)$
composé des formes paraboliques de caractère $\varepsilon$. 
Notons que par "$GAGA$ rigide analytique", ce $\HH$-module est canoniquement isomorphe
à son analogue sur $X_1(M)^{rig}/K$. La formation du $\HH$-module $S_k(M,\varepsilon,K)$
commute à l'extension des scalaires sur $K$ et la donné de $\iota$ identifie 
$S_k(M,\varepsilon,\CC_p)$ et $S_k(M,\varepsilon)$.
\par \vspace{2  mm}
$S^{D}_k(M,\varepsilon,K)$ le $K$-espace vectoriel des fonctions $f: D^*(\QQ)\backslash D^*(\AAA_f) \rightarrow
\textrm{Sym}^{k-2}(K^2)$ satisfaisant
$$f(xu)=\varepsilon(u)^{-1}u_p^{-1}f(x), \, \,
\, \, \forall (x,u) \in D^*(\AAA_f)\times U_0(M)$$
c'est un $\HH$-module de manière naturelle. La formation du $\HH$-module $S_k^D(M,\varepsilon,K)$ 
commute à l'extension des scalaires sur $K$ et la donnée de $\iota$ identifie
$S_k^D(M,\varepsilon,\CC_p)$ et $S_k^D(M,\varepsilon)$, on rappelle comment
dans ce qui suit.

\par \vspace{2 mm}A une fonction $f \in S_k^{D}(M,\varepsilon)$ est associé par définition
un morphisme $D^*(\RR)$-équivariant $\varphi_f$ de $\textrm{Sym}^{k-2}(\CC^2)^*$ vers
l'espace des fonctions complexes sur $D^*(\QQ)\backslash
D^*(\AAA)/U_1(M)$. Si $x \in D^*(\AAA)$, $$v \mapsto
\varphi_f(v)(x)$$ définit par dualité un unique élément $F_f(x) \in
\textrm{Sym}^{k-2}(\CC^2)$. On considère alors l'application qui à $f$ associe
l'élément de $S_k^{D}(M,\varepsilon,\CC_p)$ défini par  $$ x_f
\mapsto x_p^{-1}\iota^{-1}(F_f(1 \times x_f))$$ c'est l'isomorphisme cherché
comme on le vérifie immédiatement. 
\par \vspace{2 mm}

%

\vspace{1 mm}

Si $S_k(M,\varepsilon,\QQ_p)^{d-new}$ désigne le sous-espace de
$S_k(M,\varepsilon,\QQ_p)$ composé des formes
nouvelles en $d$ de $S_k(M,\varepsilon,\QQ_p)$, il vient

\begin{prop} Si $k>2$ ou $\varepsilon \neq 1$, les $\HH\otimes_{\ZZ}\QQ_p$-modules $S_k(M,\varepsilon,\QQ_p)^{d-new}$ et $S_k^{D}(M,\varepsilon,\QQ_p)$ sont
isomorphes. $S_2(M,1,\QQ_p)$ est $\HH$-isomorphe au quotient de
$S_2^{D}(M,1,\QQ_p)$ par la droite des fonctions constantes.
\end{prop}

\subsection{La correspondance $p$-adique à poids-caractère fixé} \label{cla}

	On reprend les notations du \S \,\ref{formes}, et on s'intéresse aux
systèmes de modules de Banach sur $\WW$, $$M^1:=\ME \, \, \, \, et \, \, \,
\, M^2:=\MQ$$ 
Comme on l'a vu, ils sont munis d'une action de l'algèbre $H:=\mathcal{H}$, et $U_i:=\rho_i(U_p) \in
\textrm{Comp}(M^i)$ commute à $\rho_i(H)$. On pose $X:=\{ \zeta(1+p)^k, k\geq
2, \zeta \in \mu_{p^{\infty}} \} \subset \WW(\CC_p)$,
il est très Zariski-dense dans $\WW$. On munit $M^1$ et $M^2$ de structures
classiques sur $X$ en posant pour $M^{i,cl}_{\zeta (1+p)^k}$,
$k\geq 3$, $\zeta^{p^m}=1$, 
\begin{itemize} 
\item le sous-espace de $F_{\zeta(1+p)^k,n}^{0,\varepsilon,d}$, avec
$n$ quelconque tel que $n \geq m-1$, des restrictions à
$X_1(Np^m,d)(v_n)$ des sections de $\omega^k$ convergentes sur tout
$X_1(Np^m,d)/\CC_p$, s'annulant à {\it toutes} les pointes de
$X_1(Np^m,d)/\CC_p$, si $i=1$. 
\item le sous-espace $S_k^D(Np^m,\varepsilon\chi\tau^{-k},\CC_p)$ de $F_{\zeta (1+p)^k,n}^{D,\varepsilon}$, $n \geq m-1$
quelconque, défini à la fin du paragraphe \ref{quat}.
\end{itemize}
\par \vspace{1 mm}
Ce sont bien des structures classiques, satisfaisant $(Cl)$ en $(1+p)^k\zeta$ dès que $\alpha
< k-1 $ par les assertions connues de classicité en pente petite devant le poids
(\cite{Col2} \S 8, \cite{Col3} 1.1,\cite{Buz2} \S 4). Il faut noter qu'une forme
modulaire de poids $k$ sur $X_1(Np^m,d)$ qui est propre pour $U_p$ de pente
strictement inférieure à $k\!-\!1$, et qui s'annule en toutes les pointes de
$Z_1(Np^m,d)(0)$, est en fait parabolique. On est alors dans les hypothèses de la proposition \ref{det} par le théorème \ref{JL}. On
en déduit le 


\begin{thm} \label{seriecaracteristique} Soit $h \in \HH$, alors
$\Fred_{F^{0,\varepsilon,d}}(hU_p)=\Fred_{F^{D,\varepsilon}}(hU_p)$
\end{thm}

{\it Remarques:}  \begin{itemize} \item Il est aisé de voir que le membre de droite de
cette égalité est en fait dans $1+T\Lambda\{\{T\}\}$, ainsi donc que le premier. \par \vspace{1 mm}
\item On aurait pu se restreindre, dans 
notre choix de l'ensemble $X$, à celui de $\{(1+p)^k, k\geq 3\}$ pour obtenir le même
résultat \ref{seriecaracteristique}. Via le théorème \ref{jlpfamille} (qui n'utilise
que le résultat de \ref{seriecaracteristique}), on peut alors
déduire de la propriété $(Cl)$ pour les formes quaternionique, et de
\cite{Col2}, une nouvelle preuve du théorème de contrôle de \cite{Col3}.

\end{itemize}

\par \vspace{2 mm}

Soit $\kappa \in \WW(\CC_p)$, on dispose de deux systèmes d'espaces de
Banach $F^{0,\varepsilon,d}_{\kappa}$ et
$F^{D,\varepsilon}_{\kappa}$, munis de représentations de $\HH$. Le corollaire \ref{corsemisimpl}
implique le  
\par \vspace{1 mm}

\begin{thm} \label{jlp}
$\mathcal{X}_{U_p}(\ME)=\mathcal{X}_{U_p}(\MQ)$
\end{thm}

	Autrement dit, l'espace des formes modulaires paraboliques surconvergentes de
pente finie, de poids-caractère $\kappa$, de niveau modéré $Nd$, nouvelles en $d$, de caractère
$\varepsilon$ et l'espace des formes modulaires quaternioniques $p$-adiques
de pente finie, poids-caractère $\kappa$, de niveau modéré $N$, de caractère
$\varepsilon$, ont même semi-simplification comme $\HH$-module.


\section{La correspondance en familles $p$-adiques}

\subsection{Une propriété d'unicité pour les variétés spectrales}	\label{unicite}

Soit $\WW/\CC_p$ un espace rigide réduit, soit $M$ la
donnée d'un système de modules de Banach orthonormalisables sur $\WW$, muni d'une
représentation d'une $\CC_p$-algèbre commutative $\rho: H \rightarrow \textrm{End}(M)$, et
d'un $U \in H$ tel que $\rho(U) \in \textrm{Comp}(M)$. A une telle donnée
$(M,U,H)$, on peut
attacher la série de Fredholm \Fred$_M$($U$) $\in 1+TA(\WW)\{\{T\}\}$, et l'hypersurface de Fredholm
associée $Z \subset \WW \times
\AAA^1$ définie par \Fred$_M$($U$)=0, munie de ses deux projections
$(pr_1,pr_2): Z \rightarrow \WW \times \AAA^1$. On sait alors que $Z$ a un
recouvrement admissible canonique $\mathcal{C}:=\mathcal{C}$(\Fred$_M(U)$),
composé des ouverts affinoides $\Omega \subset Z$ finis et plats sur leur image
$pr_1(\Omega)$, et ouverts fermés dans $pr_1^{-1}(pr_1(\Omega))$
(\cite{BMF} A5.8, ce qui suffit pour nos applications, \cite{Buz1} \S 4 pour
le cas général). On sait de plus construire par recollement à
l'aide de $\mathcal{C}$, un espace rigide
$D$ (\cite{eigen} \S 7 \footnote{Il s'agit de la construction "D" de la
courbe de Hecke dans \cite{eigen}, uniquement basée sur de la théorie spectrale, nous ne nous
occupons pas dans ce texte de la "C"} puis \cite{ch} \S 6), 
la "variété spectrale attachée à $(M,U,H)$", muni d'un morphisme fini
$\pi: D \rightarrow Z$. On dispose de plus d'un morphisme d'anneaux $a: H \rightarrow
A(D)$, ainsi que d'un diagramme commutatif: 


$$\xymatrix{  D  \ar@{->}^{\kappa}[dd] \ar@{->}[dr]_{\pi}
\ar@{->}[drr]^{a(U)^{-1}} 
\\  & Z \ar@{->}_{pr_2}[r] \ar@{->}^{pr_1}[dl] & \AAA_{rig}^1
\\ \WW }$$


$H$, $U \in H$ et $\WW$ étant fixés, on note $\mathcal{E}$ la catégorie dont les objets 
sont les couples $(\pi,a)$ formés d'un morphisme d'espaces rigides $\pi: D \rightarrow Z$ au dessus de
$\WW$, ainsi que d'un morphisme d'anneaux $a: H \rightarrow A(D)$ tel que
$a(U)^{-1}=pr_2\cdot \pi$. Les morphismes $\textrm{Hom}((\pi_1,a_1),(\pi_2,a_2))$ sont ceux
$(\varphi_D, \varphi_Z): \pi_1 \rightarrow \pi_2$ au dessus de $\WW$ tels que
$\forall h \in H, a_2(h) \cdot \varphi_D= a_1(h)$.
Si $X=(\pi,a)$ est un objet de $\mathcal{E}$, on notera $D(X)$ (resp. $Z(X)$) l'espace rigide $D$ source
(resp. $Z$ but) de $\pi$, $a(X):=a$, $\pi(X):=\pi$. Si
$(M,U,H)$ est comme plus haut, on notera $\mathcal{E}(M)$ l'objet de $\mathcal{E}$
associé à $M$ par la construction précédente. On notera de plus
$\mathcal{E}(M)^{\rr}$, la réduction de $\mathcal{E}(M)$, l'objet $(\pi^{\rr},a^{\rr})$ de $\mathcal{E}$, défini par  $$ \pi^{\rr}: D^{\rr}
\overset{can}{\hookrightarrow} D \overset{pi}{\longrightarrow} Z , \, \, \, 
a^{\rr}: H \overset{a}{\longrightarrow}A(D)
\overset{can}{\longrightarrow} A(D)/\textrm{Nilrad}(A(D))$$

\begin{prop} \label{spectral}Soient $(M^1,U,H)$ et $(M^2,U,H)$ comme plus haut, on suppose
$$\forall h \in H, \Fred_{M^1}(hU)=\Fred_{M^2}(hU)$$ Alors $\mathcal{E}(M^1)^{\rr}$ et
$\mathcal{E}(M^2)^{\rr}$ sont canoniquement isomorphes.
\end{prop}

{\it Preuve:} Par hypothèse, $Z(\mathcal{E}(M^1))=Z(\mathcal{E}(M^2))$, on
la note $Z$, elle est munie de sa première projection $pr_1: Z \rightarrow \WW$, on pose de plus 
$D_i=D_i(\mathcal{E}(M^i))$, $\mathcal{E}(M^i)=(\pi_i,a_i)$. 
Soit $\mathcal{C}$ le recouvrement canonique de $Z \rightarrow \WW$, $\Omega \in
\mathcal{C}$, $V:=pr_1(\Omega)$ est un ouvert affinoide de $\WW$, 
et par hypothèse $\Omega \rightarrow pr_1^{-1}(V)$ est un ouvert fermé fini sur $V$. 
Pour $n$ assez grand, $M^i(V,n)$ contient alors un sous-$A(V)$-module localement libre de rang fini,
"indépendant de
$n$", dont on note $M^i(\Omega)$ l'image dans $M_i(V)^{\dagger}$.
$M^i(\Omega)$ hérite d'une action de $H$, $\rho_i(U)$ ayant
pour polynôme caractéristique le polynôme associé à la donnée de $\Omega$
(cf. par exemple \cite{ch} 6.3.3). Soit $H_V:=H \otimes_{\CC_p} A(V)$, par construction $D_i(\Omega)$ est l'affinoide d'algèbre l'image de $H_V$
dans $\End_{A(V)}(M^i(\Omega))$. \par 
	Soit $h \in H$, montrons que si $P_{i,h}(X):=\det(h_{|M^i(\Omega)}-X.1)$, alors $P_{1,h}(X)=P_{2,h}(X)$. $V \subset
\WW$ étant réduit, il suffit de le faire après évaluation en tout $x \in
V(\CC_p)$. Soit $x \in V(\CC_p)$, alors le corollaire \ref{corsemisimpl}
s'applique à $(H,U)$ agissant sur les systèmes de modules de Banach $M^i_x$, on en conclut que les même
caractères de $H$ apparaissent dans les semi-simplifications de ces deux espaces
(avec multiplicités). $M^i(\Omega)_x$ est par définition le sous-espace de
$\bigcap_{n\in \NN} M^i_{x,n}$ sur lequel $\rho_i(U)$ a ses valeurs propres d'inverse dans
$pr_1^{-1}(\{x\})\cap \Omega$, et ces dernières ne dépendent pas de $i$, la
série caractéristique de $U$ n'en dépendant pas. Il vient donc
$P_{1,h}(X)(x)=P_{2,h}(X)(x)$, puis $P_{1,h}(X)=P_{2,h}(X)$. Si plus
généralement $h \in
H_V$, on a encore $P_{1,h}(X)=P_{2,h}(X)$ car par l'argument précédent on a
cette égalité après évaluation en tout $x \in V(\CC_p)$. \par

Soit $I_i$ l'idéal de $H_V$ noyau du morphisme $H_V \rightarrow \End_{A(V)}(M^i(\Omega))$, prouvons que $\sqrt{I_1}$=$\sqrt{I_2}$. Soit $f \in I_1$, alors
$P_{1,h}(X)=X^d$, $d$ étant le rang de $M^1(\Omega)$ (égal au rang de
$M^2(\Omega)$), on en déduit que $P_{2,h}(X)=X^d$ puis que $I_1^d \subset
I_2$. Il vient $\sqrt{I_1} \subset \sqrt{I_2}$,
puis par symétrie $\sqrt{I_1}=\sqrt{I_2}$, ce que l'on voulait. On en déduit l'existence
d'un isomorphisme d'anneau $H_V$-linéaire: $\varphi^*(\Omega):
A(D_2(\Omega)^{\rr}) \rightarrow A(D_1(\Omega)^{\rr})$. Un tel morphisme $H_V$-linéaire
est nécessairement unique s'il existe, il est au dessus de $A(\Omega)$. \par
	Soit alors $\varphi(\Omega) : D_1(\Omega)^{\rr} \rightarrow
D_2(\Omega)^{\rr}$ l'isomorphisme induit au dessus de $\Omega$. Vérifions
que si $\Omega' \subset \Omega \in \mathcal{C}$, $\varphi_{\Omega}$ envoie
$D_1(\Omega')^{\rr}$ dans $D_2(\Omega')^{\rr}$. Soit $V=pr_1(\Omega)
\subset \WW$, $V'=pr_1(\Omega')$, alors $\Omega_{V'}:=\Omega
\cap pr_1^{-1}(V') \in \mathcal{C}$ et $D_i(\Omega_{V'})$ est l'ouvert $D_i(\Omega)^{\rr} \times_V V'$
de $D_i(\Omega)^{\rr}$. $\varphi(\Omega)$ induit un $H_{V'}$-isomorphisme
$D_1(\Omega_{V'})=D_1(\Omega)^{\rr}\times_V V' \rightarrow
D_2(\Omega)^{\rr}=D_2(\Omega')^{\rr}$, ce qui conclut le cas
$\Omega'=\Omega_{V'}$. Il reste donc le cas $V=V'$. On a 
$D_i(\Omega)^{\rr}=D_i(\Omega')^{\rr}\coprod D_i(\Omega \backslash
\Omega')^{\rr}$. $\varphi_{\Omega}$ étant $A(\Omega)$ et $H_V$-linéaire sur les
fonctions, elle envoie $D_1(\Omega')^{\rr}$ isomorphiquement sur $D_2(\Omega')^{\rr}$ au dessus de $H_V$
et $A(\Omega')$. On conclut par unicité d'un tel morphisme que les
$\varphi_{\Omega}$ se recollent en un isomorphisme $D_1^{\rr} \rightarrow
D_2^{\rr}$ au dessus de $Z$ par \cite{BGR} 9.3.3/1.
%
%
$\square$ \par \vspace{2 mm} \par \vspace{2 mm}

{\it Remarque:} Il est bien sur faux en général que sous les hypothèses de
la proposition \ref{spectral}, $M^1$ et $M^2$ soient des $H$-modules
isomorphes; \ref{spectral} est la généralisation naturelle de
\ref{semisimpl}.\par \vspace{4 mm}

Nous pouvons énoncer un résultat général combinant \ref{det},
\ref{spectral} et \ref{corsemisimpl}. Fixons
$\WW$ réduit de dimension $>0$ et relativement factoriel, $(M^1,H,U)$ et
$(M^2,H,U)$ des systèmes de $H$-modules de Banach sur $\WW$
munis de structures classiques sur un sous-ensemble très Zariski-dense $X
\subset \WW(\CC_p)$ comme en \ref{critere}, satisfaisant (Cl):

\begin{thm}\label{general} Supposons pour tout $h \in H$, $x \in X$,
$$\det(1-ThU_{|M_x^{1,cl}})=\det(1-ThU_{|M_x^{2,cl}}) \in \CC_p[T]$$
alors, \par \vspace{1 mm}
\begin{itemize}
\item $\Fred_{M^1}(hU)=\Fred_{M^2}(hU) \in 1+TA(\WW)\{\{T\}\}$, \par
\vspace{1 mm}
\item $\mathcal{E}(M^1)$ est canoniquement isomorphe à $\mathcal{E}(M^2)$,
\par \vspace{1 mm}
\item Pour tout $x \in \WW(\CC_p)$, $\mathcal{X}_U(M^1_x)=\mathcal{X}_U(M^2_x)$.
\end{itemize}
\end{thm}

\par \vspace{2 mm}

	Nous allons énoncer, pour terminer ce paragraphe, un critère sur $(M,U,H)$
assurant que $D(\mathcal{E}(M))$ est réduit. Ce passage peut être omis en
première lecture, et apparaît ici par manque de référence satisfaisante. On fera les hypothèses
suivantes sur $\WW$: $\WW$ est relativement factoriel (\cite{Con} \S 4), de dimension $>0$,
et pour tout $x \in \WW(\CC_p)$, $\widehat{\OO_{X,x}}$ est intègre. On
suppose de plus que l'on dispose de $X \subset \WW(\CC_p)$ un ensemble très Zariski dense, et une structure classique sur
$X$ au sens de \S \ref{critere}, satisfaisant $(Cl)$. On fait de plus 
l'hypothèse de type "multiplicité $1$" suivante: $$ \textrm{ "Soit $\alpha \in \RR$, pour
presque tout $x \in X$, $H$ agit de manière semi-simple sur $M_x^{cl}\cap
M_x^{\leq \alpha}$}" $$

\begin{prop} \label{reduit} Sous ces hypothèses, $D(\mathcal{E}(M))$ est réduit.
\end{prop}

{\it Rappels:} Si $A$ est une algèbre affinoide, $x \in
\textrm{Specmax}(A)$, on notera $A_x$ (resp. $A^{rig}_x$) le localisé
Zariski (resp. rigide) en $x$ (\cite{BGR} 7.3.2). Ils sont tous deux
locaux noethériens, ont des complétés canoniquement isomorphes, et $A_x, A_x^{rig}$ et
$\widehat{A_x}$ sont simultanément réduits (\cite{BGR} 7.3.2/8). Notons que
si $A \rightarrow A'$ est un morphisme plat entre anneaux noethériens, alors
si $A'$ n'a pas d'idéaux premiers associés immergés, il en va de même pour
$A$. Ceci vaut en particulier pour $A_x \rightarrow A_x^{rig} \rightarrow
\widehat{A_x^{rig}}=\widehat{A_x}$. Par exemple, si tous les
$A^{rig}_x$ sont sans composantes associées immergées, alors $A$ est sans
composante associée immergée. Enfin, pour qu'un anneau noethérien $A$ sans composante associée immergée soit réduit, il
suffit que pour un ensemble $\{x_1,...,x_n\} \subset \textrm{Specmax}(A)$ tel que
chaque composante irréductible de $\textrm{Spec}(A)$ contienne un des $x_i$, chacun des
$A_{x_i}$ soit réduit. En effet, sous ces hypothèses l'application canonique $A
\rightarrow \oplus_{i=1}^r A_{x_i}$ est injective. En particulier, si un tel anneau $A$ a
son spectre irréductible, soit il est réduit, soit aucun des $A_x$, $x \in
\textrm{Specmax}(A)$ n'est réduit. 
\par 
\vspace{2 mm}

{\it Preuve:} On reprend les notations du début du \S \ref{unicite}, un
point $z \in D(\CC_p)$ sera dit classique si $\kappa(z) \in X$, et si le
caractère de $H$ obtenu par évaluation en $z$ sur $D(\CC_p)$ est dans la
semi-simplification du $H$-module $M_{\kappa(z)}^{cl}$. Les points
classiques sont alors très Zariski denses car $\WW$ est relativement factoriel
de dimension $>0$, et par (Cl). \par \vspace{1 mm}
	Soit $\Omega \in \mathcal{C}$ tel que $D(\Omega)$ contienne un point
classique, il en contient alors un ensemble Zariski-dense. Soit
$M(\Omega)$ le $A(V)$-module projectif de type fini associé à $\Omega$ comme
dans la preuve de \ref{spectral}. Soit $u \in A(D(V))\subset
\End_{A(V)}(M(\Omega))$ nilpotent, par hypothèse de multiplicité $1$, les
évaluations de $u$ sont nulles en presque tout $x$ de $X\cap V(\CC_p)$, et $u$ est donc 
nul; $D(V)$ est donc réduit s'il contient un point classique.\par \vspace{1
mm}
	Soit $\Omega \in \mathcal{C}$ quelconque, montrons que les
complétés des anneaux locaux aux points fermés de $D(\Omega)$ n'ont pas de
composantes associées immergées. C'est en fait une propriété générale 
des sous-$A$-algèbres $B$ de $\End_A(P)$ où $A$ est intègre noethérien tel que
les $\widehat{A_m}$, $m \in \textrm{Max}(A)$, soient intègres, et $P$
projectif de type fini sur $A$. En effet, si $m \in \textrm{Max}(A)$ est
fixé, $A':=\widehat{A_m}$, la platitude de $A \rightarrow A'$ entraîne que $B_m:=B
\otimes_A A'$ est canoniquement isomorphe à son image dans
$\End_{A'}(P\otimes_A A')\simeq M_r(A')$,
où $r:=\textrm{rg}_A(P)$. $A'$ étant hensélien, $B_m$ est un
produit d'algèbres locales $B_m^i$ finies sur $A'$, sans $A'$-torsion, car incluse dans
$M_r(A')$. En particulier si $Q$ est un idéal premier associé de $B_m^i$, $Q\cap A'$ est un idéal premier de $A'$ de
même hauteur que $Q$ annulant un élément de $B_m^i$, donc $Q\cap A'=0$ et
cette hauteur commune est nulle, ce que l'on voulait. \par \vspace{1 mm}
	Les $D(\Omega)$ recouvrant $D$ de manière admissible, tous les
$\OO_{D,x}^{rig}$, $x \in D$, sont sans composantes immergées associées, on conclut par le lemme suivant. $\square$ 

\begin{lemme} \label{reduc} Soit $X/\CC_p$ un espace rigide dont les anneaux
locaux n'ont pas de composantes associées immergées,
et ayant un ensemble Zariski dense de $x$ tels que $\OO_{X,x}^{rig}$ soit  
réduit, alors $X$ est réduit.
\end{lemme}

{\it Preuve:} 
Soit $X^0$ l'ensemble des $x$ de $X$ qui sont dans une seule
composante irréductible de $X$, $X^0$ est un ouvert admissible de $X$.
Soit $\red(X):=\{x \in X, \OO_{X,x}^{rig} \textrm{  est réduit
}\}$, c'est (sans condition sur $X$) un ouvert admissible de $X$, on définit
de même $\red(X^0)$. $\red(X)$ est Zariski-dense dans $X$, comme il est
ouvert il y est aussi très Zariski-dense. On en déduit que $\red(X^0)$ est
Zariski-dense dans $X^0$. Les rappels ci-dessus assurent que $\red(X^0)$ est un ouvert fermé admissible
de $X^0$, on a donc $X^0=\red(X^0)$. Si $V$ est un ouvert affinoide de $X$, $V\cap X^0$ est Zariski-dense dans $V$,
et les rappels montrent que $V \subset \red(X)$. $\square$. \par
\vspace{4 mm}

\par \vspace{2 mm}
\subsection{L'isomorphisme rigide analytique \JL$_p$}

Soient $\WW,N,p,d$ comme en \S \ref{formes},
$${D}^{0,\varepsilon,d}:=D(\mathcal{E}(F^{0,\varepsilon,d}))^{\rr}, \, \, \,
\,
{D}^{D,\varepsilon}:=D(\mathcal{E}(F^{D,\varepsilon}))^{\rr}$$

On note $Z \subset \WW \times \AAA^1$ l'hypersurface de Fredholm
associée à
\Fred$_{F^{0,\varepsilon,d}}(U_p)$=
\Fred$_{F^{D,\varepsilon}}(U_p)$ (théorème \ref{seriecaracteristique}), $a:=a(\mathcal{E}(F^{0,\varepsilon,d}))$,
$a_D:=a(\mathcal{E}(F^{D,\varepsilon}))$, on dispose aussi de
morphismes naturels $\DE \rightarrow \WW$ et $\DQ \rightarrow \WW$ que l'on
notera par le même nom $\bf{\kappa}$. 

\begin{thm} \label{jlpfamille} Il existe un unique isomorphisme rigide analytique $\JL_p$:
$D^{D,\varepsilon} \rightarrow D^{0,\varepsilon,d}$ au-dessus de $\WW$, coïncidant avec la correspondance de
Jacquet-Langlands sur les points classiques différents du point spécial. Il satisfait $\forall h \in
\HH, \, \, a(h).\JL_p=a_D(h)$.
\end{thm}

	Avant de prouver ce théorème, rappelons qu'un point $x$ de
$D^{0,\varepsilon,d}(\CC_p)$ (resp. $D^{D,\varepsilon}(\CC_p)$) est dit {\bf
classique} si : \par \vspace{1 mm}
i) $\kappa(x)=(1+p)^k\zeta$, où $k \geq 2$ est un entier,
$\zeta \in \mu_{p^{\infty}}$, \par \vspace{1 mm}

ii) le caractère $\HH \rightarrow \CC_p$ obtenu par
l'évaluation en $x$ apparaît dans la semi-simplification du $\HH$-module
$F_{(1+p)^k\zeta}^{0,\varepsilon,d,cl}$ (resp.
$F_{(1+p)^k\zeta}^{D,\varepsilon,cl}$). \par \vspace{1 mm}

	On notera $x_0$ le point de $D^{D,1}(\CC_p)$ correspondant
au caractère de $\HH$ sur la droite des fonctions constantes dans
$F_{(1+p)^2}^{D,1,cl}$, on l'appellera le {\bf point spécial}.

\par \vspace{2 mm}

{\it Preuve:} L'assertion d'unicité de \ref{jlpfamille} découle de la
Zariski-densité des points classiques de $D^{D,\varepsilon}$. Le th\'eor\`eme 
\ref{seriecaracteristique} combin\'e \`a la
proposition \ref{spectral} assure l'existence d'un isomorphisme canonique 
$\phi: \mathcal{E}(F^{D,\varepsilon})^{\rr} \rightarrow 
\mathcal{E}(F^{0,\varepsilon,d})^{\rr}$. On définit alors
$$\JL_p: \DQ \rightarrow \DE $$  
comme étant l'isomorphisme induit par $\phi$. Par construction, il est au dessus de $\WW$ et
satisfait $a=\JL_p.a_D$. \par \vspace{1 mm}

	Prouvons que $\JL_p$ induit la correspondance de Jacquet-Langlands
sur les points classiques non spéciaux. Rappelons (cf. par exemple \cite{ch} 6.2.4,
6.2.5) que si $w \in \WW(\CC_p)$, l'application qui a un point $x
\in \DE(\CC_p)$ (resp. $\DQ(\CC_p)$) tel que ${\bf \kappa}(x)=w$ associe le caractère
$\CC_p$-valué de $\HH$ d'évaluation en $x$ induit une bijection entre ${ \bf \kappa}^{-1}(w)$
et $|\mathcal{X}_{U_p}(F_{w}^{0,\varepsilon,d})|$ (resp. 
$|\mathcal{X}_{U_p}(F_{w}^{D\varepsilon})|$, voir \S \ref{semi} pour
la notation $|.|$). En particulier, si $x \in \DQ(\CC_p)$ est un point classique
non spécial, la relation $a=\JL_p.a_D$
assure que $\JL_p(x)$ correspond au même caractère de $\HH$ que $x$. Or la
correspondance de Jacquet-Langlands usuelle nous assure l'existence d'une
forme classique de poids $\kappa(x)$ ayant ce caractère sous $\HH$; par
unicité (i.e par la bijection rappelée ci-dessus) le point de $\DE(\CC_p)$
correspondant est nécessairement $\JL_p(x)$, ce qui prouve le théorème.
$\square$

\par \vspace{2 mm}

{\it Remarques:} i) Par des techniques usuelles, on peut montrer que l'adhérence $\overline{\HH}$ de la 
$\Lambda$-algèbre engendrée par l'image de $\HH$ dans $A(D^{D,\varepsilon})$ est
compacte. De ceci et de l'existence des représentations galoisiennes attachées aux formes modulaires
classiques résulte facilement (voir par exemple \cite{ch} \S 7) l'existence d'un unique
pseudo-caractère continu de dimension $2$ $$T: \textrm{Gal}(\overline{\QQ}/\QQ)_{Npd}
\longrightarrow \overline{\HH} \subset A(D^{0,\varepsilon,d})$$
tel que si $l$ est premier avec $(l,Npd)=1$, $T(\textrm{Frob}_l)=a(T_l)$. 
\par \vspace{2 mm}
ii) $\JL_p(x_0)$ correspond à la forme modulaire parabolique 
surconvergente de poids $2$, nouvelle en $d$, de $q$-développement $q+\sum_{n\geq 2}a_n q^n$
avec $a_l=l+1$ si $(l,pd)=1$, $1$ si $l|d$, et $p$ si $l=p$. Elle n'est pas
classique au sens strict précédent, mais elle est quand même convergente sur
tout $X_1(Np,d)$, bien que non cuspidale. Elle est "critique", car de pente $1$
et de poids $2$. \par \vspace{4 mm}

En prenant $X:=\{(1+p)^k, k\geq 2 \in \WW(\CC_p)\}$, les structures classiques sur
$X$ ci-dessus, $N=1$, l'hypothèse de multiplicité $1$ pour $\HH$ agissant
sur $M_{(1+p)^k}^{cl}\cap M_{(1+p)^k}^{\leq \alpha}$ sont vérifiées dès que $\frac{k-1}{2}
> \alpha $ et \ref{reduit} donne la 
\begin{prop} Si $N=1$, $D(\mathcal{E}(F^{0,\varepsilon,d}))$ et 
$D(\mathcal{E}(F^{D,\varepsilon}))
$ sont réduits. \end{prop}


\section{Quelques conséquences, remarques et questions}

\subsection{Opérateurs thêta} Soit $\kappa=(k,\chi)$ de conducteur $m$, 
$k \geq 2$, il existe un opérateur
d'entrelacements surjectif $$\C_{(1+p)^{-2}\kappa}^{\dagger} \longrightarrow
(\C_{(1+p)^{-2}{\kappa^*}}^{\dagger})\otimes {\det}^{k-1}, \, \, \, f \mapsto (\frac{d}{dT})^{k-1}(f), \, \, \, \,
\kappa^*:=(2-k,\chi)$$ en tant que représentations de $\TTT$ (cf. \cite{Buz2} \S 6, voir aussi
\cite{ST} 5.5). Son noyau est l'espace des fonctions localement polynomiales
de degré $\leq k-2$. Il induit un opérateur $\Theta^{1-k}:
F^{D,\varepsilon}_{\kappa,0}[m] \rightarrow F^{D,\varepsilon}_{\kappa,0}[m]$
avec les notations de \ref{quat}, tel que
$\forall n \in \NN, \, \, \, \,
T_n(\Theta^{1-k}(f))=n^{1-k}\Theta^{1-k}(f)$. Le noyau de $\Theta^{1-k}$ est le $\HH$-module $S_k^D(Np^m,\varepsilon\chi\tau^{-k},\CC_p)$. 
\par \vspace{1 mm}
	Du côté des courbes modulaires, on peut définir via l'application de Kodaira-Spencer
(\cite{Col2} \S 4, \cite{Col3}) un opérateur $M_{2-k}(Np^m,d)^{\dagger} \rightarrow
M_{k}(Np^m,d)^{\dagger}$ noté $\Theta^{k-1}$, qui agit sur les $q$-développements par
$(\frac{d}{dq})^{k-1}$, et satisfait donc $\forall n \in \NN, \, \,
T_n(\Theta^{k-1}(f))=n^{k-1}\Theta^{k-1}(T_n(f))$.
Les formes modulaires classiques propres, de pente finie, de poids $k$ et de
niveau $Np^md$ ne sont pas dans l'image de $\theta^{1-k}$. \par \vspace{1
mm}
	Si $x \in D^{D,\varepsilon}(\CC_p)$ est non classique mais de poids
$\kappa$, il existe une forme $f_x \in F_{\kappa}^{D,\varepsilon}[m]$
propre, non classique, ayant pour système de valeurs propres celui associé à $x$.
$\Theta^{1-k}(f_x) \in F_{\kappa^*}^{D,\varepsilon}[m]$ est propre et non
nul, et correspond donc à un point que l'on notera $\Theta^{1-k}(x) \in
D^{D,\varepsilon}(\CC_p)$. On définit de même $\Theta^{k-1}(x)$ pour tout
$x \in D^{0,\varepsilon,d}(\CC_p)$ de poids $\kappa^*$, \ref{jlpfamille}
implique la:

\begin{prop} Soit $x \in D^{D,\varepsilon}(\CC_p)$ non classique de poids
$\kappa(x)=(k,\chi)$, alors $$\Theta^{k-1}(\JL_p(\Theta^{1-k}(x)))=x$$
\end{prop}

\subsection{Questions} ${}^{}$ \par \vspace{1 mm}
(Q1) Nous n'avons pas montré que les systèmes de modules de Banach
$F^{0,\varepsilon,d}$ et $F^{D,\varepsilon}$ sont isomorphes. Est-il vrai, par
exemple, que si $\kappa \in \WW(\CC_p)$, les sous-espaces de
$(F_{\kappa}^{D,\varepsilon})^{\dagger}$ et de
$(F_{\kappa}^{0,\varepsilon,d})^{\dagger}$ composés des vecteurs de pente
finie sont des $\HH$-modules isomorphes ? Cette version "non
semi-simplifiée" de notre correspondance serait par exemple intéressante 
pour des questions de multiplicité $1$ mieux comprises du côté $GL_2$
(essentiellement par la présence du $q$-développement). \par \vspace{2 mm}

(Q2) Existe-t'il une réalisation géométrique de la correspondance donnée ?
Nous espérons revenir sur ce point dans un travail ultérieur.
\par \vspace{2 mm}

(Q3) Plusieurs autres espaces de "formes modulaires $p$-adiques" naturels peuvent être définis autant au niveau quaternionique
que pour $GL_2$. Par exemple, on peut remplacer les séries principales de
$I$ dans notre définition des formes quaternioniques $p$-adiques par des
restrictions des séries cuspidales de $GL_2(\QQ_p)$, qui apparaissent
aussi en familles sur $\WW$ et contiennent les représentations de dimension
finie usuelles aux poids-caractères arithmétiques: à quoi
correspondent-elles du côté $GL_2$ ? \par \vspace{2 mm}

Nous espérons revenir aussi sur l'étude d'une autre famille d'espaces de Banach
(construite du côté elliptique cette fois-ci), obtenue en considérant les
sections de $\omega^k$ sur la réunion finie des disques
supersinguliers de $X_1(N)/\CC_p$ ($N$ premier à $p$ disons). Ces espaces
sont aussi liés à des espaces de formes modulaires quaternioniques 
pour $D$ ramifiée en $p$ et l'infini cette fois-ci (cf. une lettre de Serre à Tate \cite{Ser2}
pour une correspondance modulo $p$).

\end{document}